\documentclass[11pt,fleqn]{article}

\usepackage[left=1.4in,right=1.4in]{geometry}

\usepackage[utf8]{inputenc}
\usepackage[T1]{fontenc}

\usepackage{amsmath,amssymb,amsthm}
\usepackage{mathtools}

\usepackage{tikz-cd}
\usetikzlibrary{decorations.pathmorphing}
\usepackage{enumerate}
\usepackage{paralist}

\usepackage{booktabs}
\usepackage{multirow}

\usepackage[hyperfootnotes=false]{hyperref}
\usepackage[capitalize,noabbrev,nameinlink]{cleveref}
\crefformat{equation}{(#2#1#3)}

\usepackage[cal=pxtx,scr=esstix,scrscaled=1.05]{mathalpha} % \mathscr font
\usepackage{aurical} % for defining \sheafhom
\usepackage{upgreek}
\usepackage{fouriernc}

\usepackage{physics}

\usepackage[backend=bibtex,style=alphabetic]{biblatex}
\DeclareFieldFormat{eprint:HAL}{\textsc{hal}: \href{https://tel.archives-ouvertes.fr/#1}{#1}}

\theoremstyle{plain}
\newtheorem{theorem}{Theorem}[section]

\newtheorem{lemma}[theorem]{Lemma}
\newtheorem{corollary}[theorem]{Corollary}

\theoremstyle{definition}
\newtheorem{definition}[theorem]{Definition}
\newtheorem{example}[theorem]{Example}

\theoremstyle{remark}
\newtheorem{remark}[theorem]{Remark}

\newcommand{\cover}{\mathcal{U}}
\newcommand{\anothercover}{\mathcal{V}}
\newcommand{\yetanothercover}{\mathcal{W}}
\newcommand{\OO}{\mathcal{O}}
\newcommand{\anotherbullet}{\star}
\newcommand{\yetanotherbullet}{\diamond}
\newcommand{\restricted}{\mathbin{\big\vert}}
\newcommand{\id}{\mathrm{id}}
\newcommand{\congto}{\xrightarrow{\raisebox{-.5ex}[0ex][0ex]{$\sim$}}}
\newcommand{\congfrom}{\xleftarrow{\raisebox{-.5ex}[0ex][0ex]{$\sim$}}}
\newcommand{\cech}{\check{\mathscr{C}}}
\newcommand{\cechd}{\hat{\mathscr{C}}}
\newcommand{\twc}{\mathfrak{a}}

\newcommand{\cst}{\mathrm{cst}}

\newcommand{\gshX}{\mathsf{Sh}(X)}
\newcommand{\gcohX}{\mathsf{Coh}(X)}
\newcommand{\gcohUX}{\mathsf{Coh}_{\cover}(X)}
\newcommand{\gcohVX}{\mathsf{Coh}_{\anothercover}(X)}
\newcommand{\gccohX}{\mathsf{CCoh}(X)}

\newcommand{\shX}{\mathsf{Sh}(\nerve{\bullet})}
\newcommand{\cartvectX}{\mathsf{Vect}^\mathrm{cart}(X_\bullet^\cover)}
\newcommand{\cartcohX}{\mathsf{Coh}^\mathrm{cart}(X_\bullet^\cover)}
\newcommand{\cartshX}{\mathsf{Sh}^\mathrm{cart}(X_\bullet^\cover)}
\newcommand{\greenX}{\mathsf{Green}(X_\bullet^\cover)}
\newcommand{\greenVX}{\mathsf{Green}(X_\bullet^\anothercover)}
\newcommand{\greenzeroX}{\mathsf{Green}_{\nabla,0}(X_\bullet^\cover)}
\newcommand{\scartvectX}{\underline{\mathsf{Vect}}^\mathrm{cart}(X_\bullet^\cover)}
\newcommand{\sgreenzeroX}{\underline{\mathsf{Green}}_{\nabla,0}(X_\bullet^\cover)}

\newcommand{\define}[1]{\textbf{#1}}
\newcommand{\nerve}[1]{X_{#1}^\cover}
\newcommand{\anothernerve}[1]{Y_{#1}^\anothercover}
\newcommand{\nervesimplex}[1]{\nerve{#1}\times\Delta^{#1}}
\newcommand{\comparison}[1]{\mathfrak{C}_{#1}}
\newcommand{\numberincircle}[1]{\text{\normalfont\Large\textcircled{\small #1}}}

\renewcommand{\leq}{\leqslant}
\renewcommand{\d}{\mathrm{d}}

\DeclareMathOperator{\End}{End}

\DeclareMathOperator{\Tot}{Tot}
\DeclareMathOperator{\HH}{H}
\DeclareMathOperator{\Ker}{Ker}
\DeclareMathOperator{\Coker}{Coker}
\DeclareMathOperator{\LL}{L}
\DeclareMathOperator{\sheafhom}{\mathscr{H}\textup{\kern -2.5pt {\large\Fontauri\slshape om}}\,}
\DeclareMathOperator{\sheafend}{\mathscr{E}\textup{\kern -2pt {\large\Fontauri\slshape nd}}\,}
\DeclareMathOperator{\footnotesizesheafend}{\mathscr{E}\textup{\kern -2pt {\Fontauri\slshape nd}}\,}
\DeclareMathOperator{\colim}{colim}
\DeclareMathOperator{\holim}{holim}
\DeclareMathOperator{\hocolim}{hocolim}

\bibliography{part-I.bib}

\title{Simplicial Chern-Weil theory for\\coherent analytic sheaves\\{\Large Part I: {\em Simplicial connections}}}
\author{Timothy Hosgood}
\date{}

\begin{document}

\maketitle

\begin{abstract}
    In ``Chern classes for coherent sheaves'', H.I.~Green constructs Chern classes in de Rham cohomology of coherent analytic sheaves.
    We construct here a formal $(\infty,1)$-categorical framework into which we can place Green's work and generalise it, also obtaining a better idea as to what exactly a \emph{simplicial connection} should be.
    The result will be the ability to work with generalised invariant polynomials (which will be introduced in the sequel to this paper) evaluated at the curvature of so-called \emph{admissible} simplicial connections to get explicit Čech representatives in de Rham cohomology of characteristic classes of coherent analytic sheaves.
\end{abstract}

\tableofcontents

\bigskip

\begin{center}
    \emph{This project has received funding from the European Research Council under the European Union’s Horizon 2020 research and innovation programme (grant agreement no.~768679), as well as ANR HODGEFUN and ANR MICROLOCAL.}
\end{center}

\medskip

{This paper is one of two to have been extracted from the author's PhD thesis \cite{Hosgood2020}. Further details and more lengthy exposition can be found there. We repeat here, however, our genuine thanks to Julien Grivaux and Damien Calaque for their tireless tutelage. We would also like to thank the reviewers for their helpful comments.}

\section{Introduction}

    \subsection{History and motivation}

    In 1980, H.I. Green, a student of O'Brian and Eells, wrote their thesis \cite{Green1980} on the subject of Chern classes of coherent sheaves on complex-analytic manifolds.
    Although the thesis was never published, an exposition was given in \cite{Toledo&Tong1986}, alongside a sketch of a proof of the Hirzebruch-Riemann-Roch formula for this construction of Chern classes.
    It combined the theory of {twisting cochains}, used with great success by Toledo, Tong, and O'Brian in multiple papers (\cite{Toledo&Tong1976,Toledo&Tong1978,OBrian&etal1981,OBrian&etal1985}), along with the {fibre integration} of Dupont (\cite{Dupont1976}), to construct, from a coherent analytic sheaf, classes in $\HH^{2k}(X,\Omega_X^{\bullet\geqslant k})$ that coincide with those given by the classical construction of Chern classes in $\HH^{2k}(X,\mathbb{Z})$ by Atiyah-Hirzebruch (\cite{Atiyah&Hirzebruch1962}).
    (It is of historical interest to mention that the approach of expressing characteristic classes in terms of transitions functions is very much in line with ideas propounded by Bott; see e.g. the subsection entitled `Concluding Remarks' in \cite[§23]{Bott&Tu1982}.)
    This construction was considered by Grivaux in his thesis \cite{Grivaux2009}, where he constructs unified Chern classes for coherent analytic sheaves (on \emph{compact} analytic manifolds) in Deligne cohomology, and where he states an axiomatisation of Chern classes that ensures uniqueness in any sufficiently nice cohomology theory (of which de Rham cohomology, as well as truncated de Rham cohomology, is an example).
    Although he states that the Grothendieck--Riemann--Roch theorem for closed immersions is not known for Green's construction of Chern classes if $X$ is non-Kähler, this turns out to not be a problem, since it follows from his other axioms by a purely formal, classical argument, involving deformation to the normal cone.

    \bigskip

    One reason that the study of Chern classes of coherent \emph{analytic} sheaves is interesting is that the techniques used in the algebraic setting seem to be entirely inapplicable to the analytic setting.
    In both the analytic and algebraic settings, Chern classes of locally free sheaves can be constructed by the splitting principle (as explained in e.g. \cite[§21]{Bott&Tu1982}) in the ``most general'' cohomology theories (Deligne--Beilinson cohomology and Chow rings, respectively); but, although coherent algebraic sheaves admit global locally free resolutions, the same is \emph{not} true of coherent analytic sheaves.
    In general, complex manifolds have very few holomorphic vector bundles, and there are whole classes of examples of coherent sheaves that do not admit a global locally free resolution (\cite[Corollary~A.5]{Voisin2002}).
    One key insight of \cite{Green1980}, however, was that the holomorphic twisting resolutions of Toledo and Tong (whose existence was guaranteed by \cite[Proposition~2.4]{Toledo&Tong1978}) could be used to construct a global resolution by ``\emph{simplicial} locally free sheaves'', or \define{locally free sheaves on the nerve}: objects that live over the Čech nerve $\nerve{\bullet}$ of a cover $\cover$ of $X$.
    The existence of such a global resolution, glued together from local pieces, is mentioned in the introduction of \cite{Hirschowitz&Simpson2001} as a problem that should be amenable to the formal theory of descent.
    Indeed, these ``simplicial sheaves'' can be constructed by taking the lax homotopy limit (in the sense of \cite[Definition~3.1]{Bergner2012}) of the diagram of model categories given by the pullback--pushforward Quillen adjunctions along the nerve of a cover of $X$.
    One very useful example of such an object is found by pulling back a global (i.e. classical) vector bundle to the nerve: given some $E\twoheadrightarrow X$, defining $\mathcal{E}^\bullet$ by $\mathcal{E}^p=(\nerve{p}\to X)^*E$.
    This actually satisfies a ``strongly cartesian'' property: it is given by the ``strict'' (i.e. not lax) homotopy limit of \cite{Bergner2012}.
    In hopefully self-explanatory notation,
    \[
        \begin{aligned}
            \operatorname*{laxholim}_{[p]\in\Delta} \mathsf{Sh}(X_p^\cover)
            &\simeq \mathsf{Sh}(X_\bullet^\cover)
        \\  \operatorname*{holim}_{[p]\in\Delta} \mathsf{Sh}(X_p^\cover)
            &\simeq \mathsf{Sh}^\mathrm{cart}(X_\bullet^\cover).
        \end{aligned}
    \]

    \bigskip

    The twisting cochains from which Green builds these resolutions are also interesting objects in their own right, having been studied extensively by Toledo, Tong, and O'Brian, as mentioned previously.
    In fact, they can be seen as specific examples of the twisted complexes of \cite{Bondal&Kapranov1991}, which are used to pretriangulate arbitrary dg-categories.
    This gives a possible moral (yet entirely informal) reason to expect the existence of resolutions such as Green's: twisted complexes give the ``smallest'' way of introducing a stable structure on a dg-category, and perfect $\OO_X$-modules can be defined as exactly the objects of the ``smallest'' stable $(\infty,1)$-category that contains $\OO_X$ and is closed under retracts.
    Alternatively, one can appeal to \cite{Wei2016}, which shows that, under certain restrictions on $(X,\OO_X)$, twisting cochains constitute a dg-enrichment of the derived category of perfect complexes.
    We also mention here one more fact, shown in \cite[Corollary~3 and Proposition~11]{BHW2017} in the language of dg-categories: analogously to how sheaves (resp. cartesian sheaves) on the nerve can be recovered as a lax (resp. not lax) homotopy limit of sheaves on each simplicial level, we can recover twisting cochains (resp. perfect twisting cochains) as a homotopy limit of locally free sheaves (resp. perfect complexes) on each simplicial level, i.e.
    \[
        \begin{aligned}
            \operatorname*{holim}_{[p]\in\Delta} \mathsf{LocFree}(X_p^\cover)
            &\simeq \mathsf{Tw}(X)
        \\  \operatorname*{holim}_{[p]\in\Delta} \mathsf{Perf}(X_p^\cover)
            &\simeq \mathsf{Tw}_\mathsf{perf}(X).
        \end{aligned}
    \]

    \bigskip

    Another problem in trying to apply Chern-Weil theory to coherent analytic sheaves is that global holomorphic (Koszul) connections rarely exist: the Atiyah class --- which coincides with the first Chern class in cohomology --- measures the obstruction of the existence of such a connection.
    The other main result of Green's thesis is the construction of simplicial connections, which are connections on ``simplicial sheaves'' pulled back along the projection $\nervesimplex{\bullet}\to\nerve{\bullet}$.
    The idea behind this construction is powerfully simple: given local connections $\nabla_\alpha$ (which always exist) on a locally-free sheaf $E$ (that is, $\nabla_\alpha$ is a connection on $E\restricted U_\alpha$), on any intersection $U_{\alpha\beta}$ we consider the path $\gamma_{\alpha\beta}(t) = t\nabla_\beta + (1-t)\nabla_\alpha$ between the two local connections as some type of connection on $\nerve{1}\times\Delta^1$.
    More generally, on $p$-fold intersections $U_{\alpha_0\ldots\alpha_p}$, we can consider the ``connection'' $\sum_{i=0}^p t_i\nabla_{\alpha_i}$ on $\nerve{p}\times\Delta^p$.
    These objects then assemble to give us what deserves to be called a \define{simplicial connection}.
    Green shows that we can take the curvature of such things, which consists of $\End(E_p)$-valued forms on $\nerve{p}\times\Delta^p$; by certain technical properties of the sheaves in his resolution, Green shows that these forms satisfy the property necessary to define a simplicial differential form (the same property as found in the equivalence relation defining the fat geometric realisation of a simplicial space), which lets us apply Dupont's fibre integration (after taking the trace, or evaluating under some other invariant polynomial) to recover (Čech representatives of) classes in de Rham cohomology.
    One thing that could be considered as missing from Green's thesis is a formal study of simplicial connections, and so this forms one of the key parts of this paper.
    It is possible to define simplicial connections in a more general setting, and study conditions that ensure that Chern-Weil theory can be applied (these give the notions of \define{admissibility}, and being \define{generated in degree zero}).
    Green's connections do indeed satisfy these formal conditions, and this provides a more rigorous reasoning for their usefulness.
    In fact, as $(\infty,1)$-categories, modulo some subtleties in the definitions, complexes of sheaves with coherent cohomology are equivalent to the homotopy colimit of so-called \define{Gre{}en complexes} endowed with simplicial connections generated in degree zero.
    This means that applying Chern-Weil theory to Gre{}en complexes does indeed give us a working version of Chern-Weil theory for complexes of sheaves with coherent cohomology.

    \subsection{Purpose and overview}

    This diptych aims to construct an $(\infty,1)$-categorical framework in which we can formally understand Green's construction, as well as showing that it agrees with other existing notions of Chern classes.
    In this first paper, we focus on the abstract definitions underlying the theory, and prove some technical results that will be used in the sequel; the next paper \cite{Hosgood2020a} will focus on explicit calculations of simplicial connections (and their curvatures) in the case of pullbacks (to the nerve) of global vector bundles.

    \medskip

    The main contributions of this paper are the following:
    \begin{enumerate}
        \item Green's construction gives a way to present a single coherent sheaf via a complex of vector bundles on the Čech nerve; in this paper we actually obtain an equivalence of $(\infty,1)$-categories between so-called ``Green complexes'' and complexes with coherent cohomology, which is then extended to an equivalence with the $(\infty,1)$-category of so-called ``Green complexes with generated-in-degree-zero simplicial connections''.
        \item The definition of a ``simplicial connection generated in degree zero'' is novel, and seems to indeed give good objects to study.
        \item Technical ideas introduced here can be applied to obtain concrete calculations (as shown in the sequel to this paper).
    \end{enumerate}

    \medskip

    In \cref{section:preliminaries} we define the necessary prerequisites (sheaves on simplicial spaces, simplicial differential forms, twisting cochains, and holomorphic twisting resolutions) and recall Green's simplicial resolution \cite[§1.4]{Green1980}, as well as exploring his example \cite[pp.~41-42]{Green1980} in more detail.
    We also define the fundamental notion of a \define{Gre{}en} complex as any complex that behaves sufficiently like some complex coming from Green's resolution.

    The purpose of \cref{section:coherent-sheaves} is to define the relevant homotopical categories, and then $(\infty,1)$-categories, of complexes of coherent sheaves (of which there are multiple possible definitions).
    We also provide (\cref{corollary:equivalence-without-connections}) the formal equivalence generalising the first part of Green's result: that complexes of sheaves with coherent cohomology are equivalent (as an $(\infty,1)$-category) to Green complexes.

    Finally, in \cref{section:simplicial-connections}, we introduce the main object of our study: \define{simplicial connections}.
    We give a geometric motivation for the definition of an \define{admissible} simplicial connection, but note that the real reason for the definition is that it is exactly the property needed in order to be able to apply a simplicial version of Chern-Weil theory (which we will do in part~II).
    We then discuss an even finer property, namely that of \define{being generated in degree zero}, and show that Gre{}en complexes always admit such simplicial connections, and that this property does indeed imply admissibility.
    This (\cref{theorem:admissible-gidz-on-Green}) is one of the two main results of this paper.
    The other main result of this paper (\cref{corollary:the-main-coherent-corollary}) is the fact that there is an equivalence of $(\infty,1)$-categories between complexes of sheaves with coherent cohomology and Green complexes endowed with generated-in-degree-zero simplicial connections, improving upon the equivalence described in \cref{section:coherent-sheaves}.
    When we develop the simplicial version of Chern-Weil theory, this equivalence will let us calculate Chern classes of coherent analytic sheaves by calculating the Chern classes of the Gre{}en complexes that resolve them.

\section{Preliminaries}\label{section:preliminaries}

    Throughout this entire paper, let $(X,\OO_X)$ be a paracompact complex-analytic manifold with its structure sheaf of holomorphic functions; let $\cover$ be a locally-finite Stein open cover of $X$ such that finite intersections are again Stein.

    \subsection{Conventions and notation}

        We write $\Delta$ to mean the abstract simplex category: its objects are the finite ordinals $[p]=[0,1,\ldots,p-1,p]$ for $p\in\mathbb{N}$; its morphisms are the order-preserving maps.
        For all $p\in\mathbb{N}$, we have, for $i\in\{0,\ldots,p-1\}$, the \define{coface} maps, which are the injections $f_p^i\colon[p-1]\to[p]$ given by omitting $i$; we also have, for $i\in\{0,\ldots,p\}$, the \define{codegeneracy} maps, which are the surjections $s_i^p\colon[p+1]\to[p]$ that send both $i$ and $i+1$ to $i$.

        We write $\Delta^\bullet$ to mean the topological simplex category: its objects are the topological simplices $\Delta^p\subset\mathbb{R}^{p+1}$ for $p\in\mathbb{N}$, where $\Delta^p$ is the set of points
        \[
            \Delta^p
            =
            \left\{
                (t_0,\ldots,t_p)
                \in \mathbb{R}^{p+1}
                \,\,\Big\vert\,\,
                t_i\geqslant0,
                \,\sum_{i=0}^p t_i =1
            \right\}
        \]
        endowed with the topology induced by the inclusion $\Delta^p\hookrightarrow\mathbb{R}^{p+1}$.
        This is an example of a \emph{co}simplicial space: it has {coface} maps $\Delta^\bullet f_p^i\colon\Delta^{p-1}\to\Delta^p$ and {codegeneracy} maps $\Delta^\bullet s_i^p\colon\Delta^{p+1}\to\Delta^p$.
        For this specific cosimplicial space, we simply write $f_p^i$ (resp. $s_i^p$) to mean $\Delta^\bullet f_p^i$ (resp. $\Delta^\bullet s_i^p$).

        Given a topological space $Y$ with open cover $\anothercover$, we define its \define{nerve} to be the simplicial space $\anothernerve{\bullet}$ given, in degree $p$, by
        \[
            \anothernerve{p} = \coprod_{\substack{\beta_0\ldots\beta_p \\ V_{\beta_0\ldots\beta_p}\neq\varnothing}} V_{\beta_0\ldots\beta_p}
        \]
        (where we write $V_{\beta_0\ldots\beta_p}$ to mean the intersection $V_{\beta_0}\cap\ldots\cap V_{\beta_p}$) and where the \define{face} maps act by
        \[
            \anothernerve{p}f_p^i\colon V_{\beta_0\ldots\beta_p} \to V_{\beta_0\ldots\widehat{\beta_i}\ldots\beta_p}
        \]
        (where the hat denotes omission) and the \define{degeneracy} maps act by
        \[
            \anothernerve{p}s_i^p\colon V_{\beta_0\ldots\beta_p} \to V_{\beta_0\ldots\beta_i\beta_i\ldots\beta_p}.
        \]

        \medskip

        We denote the \define{Čech complex} (with respect to some cover $\cover$) of an object $A$ by either $\cech_\cover^\bullet(A)$ or $\cech^\bullet(\cover,A)$.

        \medskip

        If an object has two gradings (for example, a cosimplicial object in the category of cochain complexes) then we denote one grading by $\bullet$ and the other by $\anotherbullet$ (for example, $\mathscr{F}^{\bullet,\anotherbullet}$).
        Generally, we try to keep the difference between gradings explicit in our notation, although this is sometimes at the cost of legibility when we have more than two gradings.

        \medskip

        We tend to use the terms ``vector bundle'' and ``locally free sheaf'' somewhat interchangeably.

        \medskip

        As an abuse of notation, we write $U_{\alpha_0\ldots\alpha_p}\in\cover$ to mean that each $U_{\alpha_i}$ is in $\cover$.

    \subsection{Sheaves on simplicial spaces}

        \begin{definition}\label{definition:sheaf-on-a-simplicial-space}
            Let $(Y_\bullet,\OO_{Y_\bullet})$ be a simplicial ringed space, so that each space $Y_p$ has structure sheaf $\OO_{Y_p}$.
            Then a \define{sheaf of $\OO_{Y_\bullet}$-modules on $Y_\bullet$} is a family $\mathcal{E}^\bullet$ of sheaves, where $\mathcal{E}^p$ is a sheaf of $\OO_{Y_p}$-modules on $Y_p$, along with, for all $\varphi\colon[p]\to[q]$ in $\Delta$, covariantly functorial (i.e. such that $\mathcal{E}^\bullet(\psi\circ\varphi) = \mathcal{E}^\bullet(\psi)\circ \mathcal{E}^\bullet(\varphi)$) morphisms
            \[
                \mathcal{E}^\bullet(\varphi) \colon (Y_\bullet\varphi)^*\mathcal{E}^p \to \mathcal{E}^q
            \]
            of sheaves of $\OO_{Y_q}$-modules, where we take the $\OO$-linear pullback
            \[
                (Y_\bullet\varphi)^*\mathcal{E}^p
                =
                (Y_\bullet\varphi)^{-1}\mathcal{E}^p
                \otimes_{(Y_\bullet\varphi)^{-1}\OO_{Y_p}}
                \OO_{Y_q}.
            \]

            A \define{morphism} between two such sheaves $\mathcal{E}^\bullet$ and $\mathcal{F}^\bullet$ is a family $\varphi^\bullet$ of morphisms of sheaves of $\OO_{Y_p}$-modules, where $\varphi^p\colon \mathcal{E}^p\to \mathcal{F}^p$, such that the square
            \[
                \begin{tikzcd}[row sep=huge]
                    (Y_\bullet\varphi)^*\mathcal{E}^p
                        \ar[r,"(Y_\bullet\varphi)^*\varphi^p"]
                        \ar[d,swap,"\mathcal{E}^\bullet\varphi"]
                &   (Y_\bullet\varphi)^*\mathcal{F}^p
                        \ar[d,"\mathcal{F}^\bullet\varphi"]
                \\  \mathcal{E}^q
                        \ar[r,swap,"\varphi^q"]
                &   \mathcal{F}^q
                \end{tikzcd}
            \]
            commutes.

            In the case where $Y_\bullet=\nerve{\bullet}$ we often say, as a mild abuse of language, \define{vector bundle on the nerve} to mean ``\emph{locally free} sheaf of $\OO_{\nerve{\bullet}}$-modules''.
        \end{definition}

        \begin{remark}\label{remark:simplicial-sheaves-terminology}
            We refrain from calling such objects simplicial sheaves, since this phrase can be justifiably interpreted in various ways.
            They are indeed simplicial objects in some category of sheaves, but this is maybe not entirely immediately obvious, since they look more cosimplicial at first glance.
            They can also be expressed as (co)lax limit objects.\footnote{We can formalise this using the lax homotopy limit (of \cite[Definition~3.1]{Bergner2012}) of the diagram of model categories given by the pullback-pushforward Quillen adjunctions along the $\nerve{\bullet}$, which, in the notation of \cite{Bergner2012}, means taking $F_{\alpha,\beta}^\theta = (\nerve{\bullet}\theta)^*$ and $u_{\alpha,\beta}^\theta = \mathcal{E}^\bullet(\theta)$.}
            For more details, we refer the interested reader to \cite{Hosgood2024}.

            We repeat the fact that these objects are \emph{covariant} with respect to $\Delta$ (that is, \emph{cosimplicial objects}), since they are contravariant objects (sheaves) on a contravariant space (a simplicial space).
        \end{remark}

        \begin{remark}\label{remark:rank-of-vector-bundle-on-the-nerve}
            There are no conditions on the ranks of a vector bundle on the nerve: it could be the case that $\mathcal{E}^p$ is of rank $\mathfrak{r}$ but $\mathcal{E}^q$ is of rank $\mathfrak{s}$.
            By definition, however, the rank \emph{is} constant over different open sets of the same simplicial degree: $\mathcal{E}^p\restricted U_{\alpha_0\ldots\alpha_p}$ is of the same rank as $\mathcal{E}^p\restricted U_{\beta_0\ldots\beta_p}$.
        \end{remark}

        \begin{definition}
            A sheaf $\mathcal{E}^\bullet$ on a simplicial space is said to be \define{strongly cartesian} if the $\mathcal{E}^\bullet\varphi$ are isomorphisms for all $\varphi\colon[p]\to[q]$; a globally bounded complex $\mathcal{E}^{\bullet,\anotherbullet}$ of sheaves on a simplicial space is said to be \define{cartesian} if all the $\mathcal{E}^\bullet\varphi$ are quasi-isomorphisms, and \define{strongly cartesian} if they are all (strict) isomorphisms.

            This former condition (of being cartesian) strengthens the ``colax limit object'' description mentioned in \cref{remark:simplicial-sheaves-terminology} to ``\emph{strict} limit object''.
        \end{definition}

        \begin{example}\label{example:global-vector-bundles}
            Given some locally free sheaf $E$ on $X$, we can pull it back to the nerve: define $E^p=(\nerve{p}\to X)^*E$.
            We call such sheaves \define{global vector bundles on the nerve}.
            This gives a particularly well-behaved family of examples: all global vector bundles on the nerve have \emph{constant rank} (across simplicial levels) and are \emph{strongly cartesian} (in fact, all the $E^\bullet\varphi$ are identity maps, which is even stronger).
        \end{example}

    \subsection{Simplicial differential forms}

        \begin{definition}
            Let $Y_\bullet$ be a simplicial complex manifold.
            Following \cite{Dupont1976}, we define a \define{simplicial differential $r$-form on $Y_{\bullet}$} to be a family $\omega_\bullet$ of forms, with $\omega_p$ a global section of the sheaf
            \[
                \bigoplus_{i+j=r}
                    \pi_{Y_p}^* \Omega_{Y_p}^i
                    \otimes_{\OO_{Y_p\times \Delta^p_\mathrm{extd}}}
                    \pi_{\Delta^p_\mathrm{extd}}^* \Omega_{\Delta^p_\mathrm{extd}}^j
            \]
            (where $\Delta^p_\mathrm{extd}$ is the affine subspace of $\mathbb{R}^{p+1}$ given by the vanishing of $1-\sum_{m=0}^p x_p$; where $\Omega_{Y_p}$ is the sheaf of \emph{holomorphic} forms, and $\Omega_{\Delta^p_\mathrm{extd}}$ is the sheaf of \emph{smooth} forms; and where $\pi_{Y_p}$ and $\pi_{\Delta^p_\mathrm{extd}}$ are the projection maps from $Y_p\times\Delta^p_\mathrm{extd}$) such that, for all \emph{coface} maps $f_p^i\colon[p-1]\to[p]$,
            \begin{equation}
            \label{equation:simplicial-gluing-condition-for-forms}
                \left(Y_\bullet f_p^i\times\id\right)^*\omega_{p-1}
                = \left(\id\times f_p^i\right)^*\omega_p
                \in \Omega^r(Y_{p}\times\Delta^{p-1}).
            \end{equation}

            We write $\Omega^{r,\Delta}(Y_\bullet)$ to mean the algebra of all simplicial differential $r$-forms on $Y_\bullet$.
            We can describe each $\omega_p$ as a form \define{of type $(i,j)$}, by writing $\omega_p = \xi_p\otimes\tau_p$, where $\xi_p$ is the $Y_p$-part of $\omega_p$, and $\tau_p$ is the $\Delta^p$-part of $\omega_p$; then $i=|\xi_p|$ and $j=|\tau_p|$.
            This lets us define a \define{differential}
            \[
                \d \colon \Omega^{r,\Delta}(Y_\bullet) \longrightarrow \Omega^{r+1,\Delta}(Y_\bullet)
            \]
            which is given by the Koszul convention with respect to the type of the form:
            \begin{align*}
                \d(\xi_p\otimes\tau_p)
                &= \left(\d_{Y_\bullet} + (-1)^{|\xi_p|}\d_{\Delta^\bullet}\right) (\xi_p\otimes\tau_p)
            \\  &= \d\xi_p\otimes\tau_p + (-1)^{|\xi_p|}\xi_p\otimes\d\tau_p.
            \end{align*}
        \end{definition}

        \begin{remark}
            The condition in \cref{equation:simplicial-gluing-condition-for-forms} can be understood as asking that the forms descend to the fat geometric realisation (as explained in \cite{Dupont1976}), or in terms of framings (as alluded to in \cite{Hosgood2020}).
        \end{remark}

        \begin{lemma}[Dupont's fibre integration]\label{lemma:dupont's-fibre-integration}
            There is a \emph{quasi-isomorphism} which, for each fixed degree $r$, consists of the map
            \begin{equation}
                \int_{\Delta^\bullet}\colon \Omega^{r,\Delta}({Y_\bullet})
                \to
                \bigoplus_{p=0}^r\Omega^{r-p}({Y_p})
            \end{equation}
            induced by \define{fibre integration}
            \begin{equation}
                \int_{\Delta^p}\colon\Omega^{r,\Delta}({Y_\bullet})
                \to
                \Omega^{r-p}({Y_p})
            \end{equation}
            where the latter is given by integrating the type-$(r-p,p)$ part of a simplicial form over the geometric realisation of the $p$-simplex with its canonical orientation.
            \begin{proof}
                The proof in the \emph{smooth} case is exactly \cite[Theorem~2.3]{Dupont1976} along with \cite[Remark~1, §2]{Dupont1976}; the proof in the holomorphic case works almost identically.
                More details (which are especially useful for explicit calculations) can be found in \cite{Hosgood2020}.
            \end{proof}
        \end{lemma}

        \begin{remark}
            There is a possibility for confusion (and many sign errors) here: \cite{Dupont1976} uses the convention of writing simplicial differential forms as forms on $\Delta^p\times Y_{p}$; \cite{Green1980} does the opposite, writing $Y_{p}\times\Delta^p$.
            We opt for the latter.

            However, this does not really concern us, given our present purposes: we do not perform any actual calculations with fibre integration in this paper; we mention it only to be able to state \cref{example:fibre-integration-gives-cech-de-rham}, which tells us that `having characteristic classes defined at the level of simplicial differential forms will let us recover classes defined at the level of de Rham cohomology'.
            This will form the basis of the sequel to this paper.
        \end{remark}

        \begin{example}\label{example:fibre-integration-gives-cech-de-rham}
            Taking $Y_\bullet=\nerve{\bullet}$ gives
            \[
                \int_{\Delta^\bullet}\colon \Omega^{r,\Delta}({\nerve{\bullet}})
                \to
                \bigoplus_{p=0}^r\Omega^{r-p}({\nerve{p}})
                \cong
                \Tot^r\cech_\cover^\bullet(\Omega^\bullet_X)
            \]
            where $\cech$ denotes the Čech complex.
            It is interesting to note that the conditions (namely being locally finite and Stein) that we impose on $\cover$ are really (ar far as fibre integration is concerned) only to ensure that this quasi-isomorphism will eventually calculate de-Rham cohomology; the \emph{existence} of the quasi-isomorphism in \cref{lemma:dupont's-fibre-integration} does \emph{not} depend on the properties of the cover.
        \end{example}

    \subsection{Green's resolution}

        \begin{definition}
            A (cochain) complex $K^\bullet$ of sheaves on a locally-ringed space $(X,\OO_X)$ is said to be \define{perfect} if, locally, it is quasi-isomorphic to a bounded complex of locally free sheaves.
            That is, if, for every point $x\in X$, there exists some open neighbourhood $U$ of $x$ and some bounded complex $E^\bullet_U$ of locally free sheaves on $U$ such that $K^\bullet\restricted U$ is quasi-isomorphic to $E^\bullet_U$.
        \end{definition}

        \begin{definition}\label{definition:deleted-cech-etc}
            Suppose that, over each $U_\alpha\in\cover$, we have a finite-length complex $(V_\alpha^\bullet,\d_\alpha)$ of locally free $\OO_{U_\alpha}$-modules, the collection of which we refer to simply as $V^\bullet$.
            Define the collection $\End^q(V)$ of \define{degree-$q$ endomorphisms of $V$} over each $U_{\alpha_0\ldots\alpha_p}$ by
            \begin{equation*}
                \End^q(V)\restricted{U_{\alpha_0\ldots\alpha_p}} = \Big\{\big(f^i\colon V_{\alpha_p}^i\restricted{U_{\alpha_0\ldots\alpha_p}}\to V_{\alpha_0}^{i+q}\restricted{U_{\alpha_0\ldots\alpha_p}}\big)_{i\in\mathbb{Z}} \,\Big|\, \d_{\alpha_p}\circ\,f^i = f^{i+1}\circ\,\d_{\alpha_0}\Big\}.
            \end{equation*}
            That is, an element of $\End^q(V)$ is a `true' (in that it commutes with the differentials) morphism of degree $q$ of chain complexes $V_{\alpha_p}^\bullet\to V_{\alpha_0}^\bullet$.

            Following \cite[0.A]{Green1980}, we define the \define{deleted Čech complex}
            \begin{equation*}
                \cechd^p\big(\cover,\End^q(V)\big)
                =
                \smashoperator{\prod_{\substack{(\alpha_0\ldots\alpha_p)\\\text{s.t. }U_{\alpha_0\ldots\alpha_p}\neq\varnothing}}}
                \End^q(V)\restricted {U_{\alpha_0\ldots\alpha_p}}
            \end{equation*}
            with \define{deleted Čech differential} $\hat{\delta}\colon\cechd^p\big(\cover,\End^q(V)\big) \to \cechd^{p+1}\big(\cover,\End^q(V)\big)$ given by
                \begin{equation*}
                    (\hat{\delta}c)_{\alpha_0\ldots\alpha_{p+1}} = \sum_{i=1}^p(-1)^ic_{\alpha_0\ldots\widehat{\alpha_i}\ldots\alpha_{p+1}}\restricted {U_{\alpha_0\ldots\alpha_{p+1}}}
                \end{equation*}
            (note that the sum is only over $i\in\{1,\ldots,p\}$, missing out both $i=0$ \textit{and} $i=p+1$).
        \end{definition}

        \begin{definition}\label{definition:holomorphic-twisted-cochain}
            A \define{holomorphic twisting cochain for $(\cover,V^\bullet)$}, where $\cover$ and $V^\bullet$ are as in \cref{definition:deleted-cech-etc}, is an element
            \begin{equation*}
                \twc=\sum_{k\geqslant0}\twc^{k,1-k}\in\Tot^1\cechd^\bullet\big(\cover,\End^\anotherbullet(V)\big)
            \end{equation*}
            such that
            \begin{enumerate}[(i)]
                \item $\twc^{1,0}_{\alpha\alpha}=\id$;
                \item $\twc^{0,1}_\alpha=\d_\alpha$;
                \item $\twc$ satisfies the Maurer-Cartan equation:
                    \[
                        \hat{\delta}\twc+\twc\cdot\twc=0.
                    \]
            \end{enumerate}
        \end{definition}

        \begin{definition}\label{definition:holomorphic-twisted-resolution}
            Let $F$ be a sheaf of $\OO_X$-modules on $X$.
            Then a \define{holomorphic twisting resolution of $F$} is a triple $(\cover,V^\bullet,\twc)$ such that the following conditions are satisfied:
            \begin{enumerate}[(i)]
                \item $\cover=\{U_\alpha\}$ is a locally-finite Stein open cover of $X$;
                \item $V^\bullet=(V^\bullet_\alpha,\d_\alpha)$ is a collection of local locally free resolutions of $F$ over each $U_\alpha$, of \emph{globally} bounded length, i.e. each $V^\bullet_\alpha$ is a resolution of $F\restricted {U_\alpha}$ by locally free $\OO_{U_\alpha}$-modules, and there exists some $B\in\mathbb{N}$ such that \emph{every} $V^\bullet_\alpha$ is of length no greater than $B$;
                \item $\twc$ is a holomorphic twisting cochain for $(\cover,V^\bullet)$ \textit{over $F$} --- that is, $\twc$ is a holomorphic twisting cochain for $(\cover,V^\bullet)$ such that we have the following commutative diagram:
                    \begin{equation*}
                        \begin{tikzcd}[column sep=small]
                            V^\bullet
                                \ar{rr}{\twc^{1,0}}
                                \ar{dr}
                        &&  V^\bullet
                                \ar{dl}
                        \\& F
                        \end{tikzcd}
                    \end{equation*}
                \item on degenerate simplices \emph{of the specific form $\alpha=(\alpha_0\ldots\alpha_p)$ with $\alpha_i=\alpha_{i+1}$ for some~$i$}, we have that $\twc^{k,1-k}=0$ for $k\geqslant2$.
            \end{enumerate}
        \end{definition}

        \begin{remark}
            There are a few existence criteria for holomorphic twisting resolutions.
            In particular, when $F$ is coherent, \cite[Lemma~8.13]{Toledo&Tong1976} and \cite[Lemma~2.4]{Toledo&Tong1978} (with the latter being an applied version of the former) both show that a holomorphic twisting cochain exists (and the latter actually shows a stronger result using the Hilbert syzygy theorem, namely that we can ensure that our global-length bound $B$ is no more than the dimension of $X$).
        \end{remark}

        \begin{definition}\label{definition:elementary-sequence}
            Given a ring $R$ and some $R$-modules $N_1,\ldots,N_s$, we say that a sequence $0\to M_r\to\ldots\to M_0\to0$ of $R$-modules is \define{$(N_1,\ldots,N_s)$-elementary} if it is a direct sum of sequences of the form
            \[
                (0\to N_i\xrightarrow{\id}N_i\to 0)[n]
            \]
            for some $i\in\{1,\ldots,s\}$ and $n\in\mathbb{Z}$.
            Given complexes $V_1^\bullet,\ldots,V_t^\bullet$ of $R$-modules, we say that a sequence is \define{$(V_1^\bullet,\ldots,V_t^\bullet)$-elementary} if it is $\mathcal{N}$-elementary, where
            \[
                \mathcal{N} = \left\{V_i^j \,\,\big|\,\, 1\leq i\leq t,1\leq j\leq s\right\}.
            \]
        \end{definition}

        \begin{theorem}[Green's resolution]\label{theorem:green's-resolution}
            Let $\mathscr{F}$ be a coherent sheaf of $\OO_X$-modules on a paracompact complex-analytic manifold $X$ with locally-finite Stein cover $\cover=\{U_\alpha\}_{\alpha\in I}$.
            Let $(\cover,V^\bullet,\twc)$ be a holomorphic twisted resolution of $\mathscr{F}$.
            Denote by $\mathscr{F}^\bullet$ the pullback of $\mathscr{F}$ to the nerve $\nerve{\bullet}$.
            Then there exists a resolution $\mathcal{E}^{\bullet,\anotherbullet}$ of $\mathscr{F}^\bullet$ (in the sense that the morphism $\mathcal{E}^{p,\anotherbullet}\to\mathscr{F}^p$ is a quasi-isomorphism in for every $p\in\mathbb{N}$) by vector bundles on the nerve:
            \begin{equation*}
                0 \to \mathcal{E}^{\bullet,0} \to \ldots \to \mathcal{E}^{\bullet,n} \to \mathscr{F}^\bullet
            \end{equation*}
            where $n=\dim X$.
            Further, each $\mathcal{E}^{\bullet,j}$ satisfies the following properties:
            \begin{enumerate}[(i)]
                \item $\mathcal{E}^{0,\anotherbullet}\restricted {U_\alpha} \cong V^\anotherbullet_\alpha$;
                \item for all coface\footnote{In fact, in \cite[§1.4]{Green1980}, properties (ii) and (iii) are stated for arbitrary compositions of coface (resp. codegeneracy) maps instead of simply for single coface (resp. codegeneracy) maps.} maps $f_p^i\colon[p-1]\to[p]$, the map
                    \[
                        \mathcal{E}^{\bullet,\anotherbullet} f_p^i\colon \big(\nerve{\bullet} f_p^i\big)^*\mathcal{E}^{p-1,\anotherbullet}\to\mathcal{E}^{p,\anotherbullet}
                    \]
                    of complexes of sheaves on $\nerve{p}$ is injective, and $\Coker\left(\mathcal{E}^{\bullet,\anotherbullet} f_p^i\right)$ is an elementary sequence;\footnote{Here the complex is indexed by $\anotherbullet$ from $0$ to $n$, and $\mathcal{E}^{\bullet,\anotherbullet} f_p^i$ is a map of complexes.}
                \item for all codegeneracy maps $s_i^p\colon[p+1]\to[p]$, the map
                    \[
                        \mathcal{E}^{\bullet,\anotherbullet} s_i^p\colon \big(\nerve{\bullet} s_i^p\big)^*\mathcal{E}^{p+1,\anotherbullet}\to\mathcal{E}^{p,\anotherbullet}
                    \]
                    of complexes of sheaves on $\nerve{p}$ is surjective, and $\Ker\left(\mathcal{E}^{\bullet,\anotherbullet} s_i^p\right)$ is an elementary sequence.
            \end{enumerate}
            These follow from the fact that, for all ${\upgamma}\leqslant{\upbeta}\leqslant{\upalpha}=(\alpha_0\ldots\alpha_p)$, writing $\mathcal{E}_{\upalpha}$ to mean $\mathcal{E}^p\restricted U_{\upalpha}$, we have the following:
            \begin{enumerate}[(i)]
                \setcounter{enumi}{3}
                \item $\mathcal{E}_{\upalpha}^\anotherbullet \cong \mathcal{E}^\anotherbullet_{\upbeta}\restricted U_{\upalpha}\oplus \mathcal{L}_{{\upalpha},{\upbeta}}^\anotherbullet$ for some $(V^\anotherbullet_{\alpha_0},\ldots,V^\anotherbullet_{\alpha_p})$-elementary sequence $\mathcal{L}_{{\upalpha},{\upbeta}}^\anotherbullet$;
                \item $\mathcal{L}_{{\upalpha},{\upgamma}}^\anotherbullet \cong \mathcal{L}_{{\upbeta},{\upgamma}}^\anotherbullet\restricted U_{\upalpha}\oplus \mathcal{L}_{{\upalpha},{\upbeta}}^\anotherbullet$;
                \item over each $U_{\upalpha}$ there is the commutative diagram
                    \begin{equation*}
                        \begin{tikzcd}
                            0 \rar
                            &\mathcal{E}_{\upbeta}^\anotherbullet \rar
                            &\mathcal{E}_{\upalpha}^\anotherbullet \rar
                            &\mathcal{L}_{{\upalpha},{\upbeta}}^\anotherbullet \rar
                            &0\\
                            0 \rar
                            &\mathcal{E}_{\upgamma}^\anotherbullet\oplus \mathcal{L}_{{\upbeta},{\upgamma}}^\anotherbullet \rar \uar{\wr}
                            &\mathcal{E}_{\upgamma}^\anotherbullet\oplus \mathcal{L}_{{\upalpha},{\upgamma}}^\anotherbullet \rar \uar{\wr}
                            &\mathcal{L}_{{\upalpha},{\upbeta}}^\anotherbullet \rar \uar{\id}
                            &0
                        \end{tikzcd}
                    \end{equation*}
                    (omitting the restriction notation), where the bottom map is induced by the natural inclusion $\mathcal{L}_{{\upbeta},{\upgamma}}^\anotherbullet\hookrightarrow \mathcal{L}_{{\upalpha},{\upgamma}}^\anotherbullet$ coming from $\mathcal{L}_{{\upalpha},{\upgamma}}^\anotherbullet \cong \mathcal{L}_{{\upbeta},{\upgamma}}^\anotherbullet\oplus \mathcal{L}_{{\upalpha},{\upbeta}}^\anotherbullet$.
                \end{enumerate}

            \begin{proof}
                This is \cite[§1.4]{Green1980}.
            \end{proof}
        \end{theorem}

        \begin{definition}\label{definition:Gre{}en}
            We say that any complex $\mathcal{E}^{\bullet,\anotherbullet}$ of vector bundles on the nerve satisfying conditions (iv), (v), and (vi) of \cref{theorem:green's-resolution} is \define{Gre{}en}.
            We say that a single vector bundle on the nerve $\mathcal{E}^\bullet$ is Gre{}en if the complex $\mathcal{E}^\bullet[0]$ concentrated in degree zero is Gre{}en.

            Note that we do not prove whether or not every Green complex actually comes from Green's resolution applied to a coherent sheaf, but this will not matter.
        \end{definition}

        \begin{remark}
            Both \cref{definition:holomorphic-twisted-resolution} and \cref{theorem:green's-resolution} can be generalised to \emph{complexes}\footnote{Since we are in the analytic setting, there is some subtlety surrounding the `good' definition of the category of complexes of coherent sheaves. We discuss this further in \cref{section:coherent-sheaves}.} of coherent sheaves: the fact that Green's resolution still works follows exactly the same lines as the original proof (although we explain the details in \cite{Hosgood2020}); the fact that holomorphic twisting resolutions of complexes of coherent sheaves exist is explained in the proof of \cref{lemma:green-gives-essential-surjectivity}.
        \end{remark}

        \begin{corollary}\label{corollary:green-gives-cartesian-things}
            Green complexes are, in particular, \emph{cartesian} complexes of locally free sheaves on the nerve.
            \begin{proof}
                Since taking the direct sum with an elementary sequence is a quasi-isomorphism, this follows from properties~(ii) and (iii) of \cref{theorem:green's-resolution}.
            \end{proof}
        \end{corollary}

        \begin{remark}
            The key point to make here is that Green's simplicial resolution is \textbf{not} just the data of a resolution, but also the properties governing the (co)kernels.
            This is important because it will tell us, in particular, that we get \emph{admissible} simplicial connections on each sheaf in the resolution, which is vital in defining characteristic classes.
            We use the coface-injectivity property when discussing admissibility of simplicial connections, but don't seem to need the codegeneracy-surjectivity property anywhere.
            This might be because the simplicial condition for e.g. simplicial differential forms relies only on (co)face maps, and so we don't care so much about (co)degeneracy maps in our formalism, and should really be talking about $\Delta$- (or semi-simplicial-) sets.
        \end{remark}

        \begin{remark}
            Green's resolution is \emph{not} functorial, but this does not overly matter to us: we only use it to show essential surjectivity of a functor (in \cref{lemma:green-gives-essential-surjectivity}), and so only care about what it does to \emph{objects}.
        \end{remark}

        \begin{example}[Green's example]\label{example:green's-example}
            Let $X=\mathbb{P}_\mathbb{C}^1\cong\mathbb{C}\cup\{\infty\}$ be the projective line over $\mathbb{C}$, with cover $\cover=\{U_1,U_2\}$, where $U_1=X\setminus\{\infty\}$ and $U_2=X\setminus\{0\}$.
            Let $\mathscr{F}$ be the coherent sheaf given by $\OO_X/\mathscr{I}$, where $\mathscr{I}=\mathbb{I}(\{0\})$ is the sheaf of ideals corresponding to the subvariety $\{0\}\subset X$.
            We use $\alpha$ and $\beta$ to denote indices in the indexing set $\{1,2\}$ with the tacit assumption that, whenever we do so, $\alpha\neq\beta$.

            The stalks of $\mathscr{F}$ are easy enough to understand:
            \[
                \mathscr{F}_z =
                \begin{cases}
                    \mathbb{C} &\mbox{if $z=0$};\\
                    0 &\mbox{if $z\neq0$}.
                \end{cases}
            \]
            This makes it easy to find local resolutions by $\OO_{U_\alpha}$-modules: over $U_1$ we have the resolution
            \begin{equation*}
                (\xi_1^\bullet\to \mathscr{F}\restricted U_1) =  \underbrace{\big(0\to\OO_X\restricted {U_1}\xrightarrow{f\mapsto z\cdot f}\OO_X\restricted {U_1}\big)}_{\xi_1^\bullet}\to \mathscr{F}\restricted {U_1},
            \end{equation*}
            and over $U_2$ we have the resolution
            \begin{equation*}
                (\xi_2^\bullet\to \mathscr{F}\restricted U_2) = \underbrace{\big(0\to0\to0\big)}_{\xi_2^\bullet}\to \mathscr{F}\restricted {U_2}.
            \end{equation*}

            Since neither $0$ nor $\infty$ are in $U_{12}$, the map $f\mapsto z\cdot f$ gives an automorphism of $\OO_X\restricted U_{12}$.
            This means that the homology of $\xi_1\restricted U_{12}$ is zero, and hence isomorphic to the homology of $\xi_2\restricted U_{12}$.
            The zero map thus gives a quasi-isomorphism $\xi_1\restricted U_{12} \congto \xi_2\restricted U_{12}$, and higher homotopies can be constructed by an inductive method.

            To construct a resolution by locally free sheaves on the nerve we need, in particular, a resolution of $\mathscr{F}\restricted U_{12}=0$ by locally free sheaves over $U_{12}$.
            Now we need to make a choice, since we have two different (but quasi-isomorphic) possibilities for such a complex:
            \begin{align*}
                \,^{(1)}\xi_{12}^\bullet &= \big(0\to\OO_X\vert{U_{12}}\xrightarrow{f\mapsto z\cdot f}\OO_X\restricted {U_{12}}\big)\to \mathscr{F}\restricted {U_{12}},\\
                \,^{(2)}\xi_{12}^\bullet &= \big(0\to0\to0\big)\to \mathscr{F}\restricted {U_{12}}.
            \end{align*}
            But we see that adding the elementary sequence
            \begin{equation*}
                \,_{(2)}^{(1)}\mathcal{L}^\bullet=(0\to\OO_X\xrightarrow{\id}\OO_X\to0)
            \end{equation*}
            to the latter gives us something isomorphic (since neither $0$ nor $\infty$ are in $U_{\alpha\beta}$, so both $1/z$ and $z$ are well defined and holomorphic) to the former:
            \begin{equation*}
                \begin{tikzcd}
                    \qquad\quad\,\,^{(1)}\xi_{12}^\bullet=
                        \big(0 \ar[r]
                        &\OO_X\restricted {U_{12}} \ar[r,"f\mapsto z\cdot f"] \ar[d,swap,"f\mapsto z\cdot f"]
                        &\OO_X\restricted {U_{12}}\big) \ar[d,swap,"\id"]\\
                    \,^{(2)}\xi_{12}^\bullet\oplus\,_{(2)}^{(1)}\mathcal{L}^\bullet =
                        \big(0 \ar[r]
                        &\OO_X\restricted {U_{12}} \ar[r,"\id"]
                        &\OO_X\restricted {U_{12}}\big).
                \end{tikzcd}
            \end{equation*}
            We denote this isomorphism by
            \begin{equation*}
                A_{12}^\bullet\colon\,^{(1)}\xi_{12}^\bullet\congto\,^{(2)}\xi_{12}^\bullet\oplus\,_{(2)}^{(1)}\mathcal{L}^\bullet
            \end{equation*}
            i.e. $A_{12}^0=\id$ and $A_{12}^1=(z\cdot-\,)$.
            We don't actually need this second complex (the codomain of $A_{12}^\bullet$) to construct the resolution, but it will prove useful in \cref{example:green's-example-cont}.

            At this point, the construction stabilises: since our cover consists of only two distinct opens sets, any $p$-intersection for $p\geqslant3$ will be exactly some $2$-intersection (either $U_{11}$, $U_{12}$, or $U_{22}$).
            We have thus constructed a complex of $\OO_{\nerve{\bullet}}$-modules on the nerve of $X$ that give a resolution of $\mathscr{F}$ (pulled back to the nerve):
            \begin{equation*}
                0\to\mathcal{E}^{\bullet,1}\to\mathcal{E}^{\bullet,0}\to\mathscr{F}^\bullet
            \end{equation*}
            where $\mathcal{E}^{\bullet,1}$ and $\mathcal{E}^{\bullet,0}$ are equal, and defined as
            \begin{align*}
                &\mathcal{E}^{0,i}=\OO_X\restricted {U_1}\sqcup0 \quad \text{over}\quad U_1\sqcup U_2\\
                &\mathcal{E}^{1,i}=\OO_X\restricted {U_{12}} \qquad \text{over}\quad U_{12}
            \end{align*}
            and the map $\mathcal{E}^{\bullet,1}\to\mathcal{E}^{\bullet,0}$ is given by
            \begin{align*}
                \big((z\cdot-\,)\sqcup0\big) \colon \mathcal{E}^{0,1}&\to\mathcal{E}^{0,0} \quad\text{over}\quad U_1\sqcup U_2\\
                (z\cdot-\,) \colon \mathcal{E}^{1,1}&\to\mathcal{E}^{1,0} \quad\text{over}\quad U_{12}.
            \end{align*}
        \end{example}

\section{Coherent sheaves}\label{section:coherent-sheaves}

    \subsection{Homotopical categories}

        \begin{definition}
            In this chapter we are interested in \define{relative categories}: pairs $(\mathcal{C},\mathcal{W})$ where $\mathcal{C}$ is a category and $\mathcal{W}$ (whose morphisms we call \emph{weak equivalences}) is a wide subcategory of $\mathcal{C}$.
            A relative category is said to be a \define{homotopical category} if its weak equivalences satisfy the {2-out-of-6 property} \cite[\S3]{Rie2020}.
            We often write a relative (or homotopical) category $(\mathcal{C},\mathcal{W})$ simply as $\mathcal{C}$, omitting the weak equivalences from our notation.

            Using the formalism of \cite{Rezk2000} along with the results of \cite{Barwick&Kan2013}, we can think of a homotopical category $(\mathcal{C},\mathcal{W})$ as presenting the $(\infty,1)$-category $\LL_\mathcal{W}\mathcal{C}$, which is the complete Segal space given by taking a Reedy fibrant replacement of the Rezk/simplicial nerve $N(\mathcal{C},\mathcal{W})$.
            In particular, \cite[§1.2~(ii)]{Barwick&Kan2013} tells us that any homotopically full relative subcategory of a partial model category (defined loc. cit.) is again a partial model category, and all of the categories that we study here are such subcategories of either the category of complexes of sheaves on $X$ or of the category of complexes of sheaves on $\nerve{\bullet}$, both of which\footnote{The former having, for example, the projective model structure from \cite{Hovey2001}; the latter having the model structure coming from its construction as a lax homotopy limit via the formalism of \cite{Bergner2012}, as described in the footnote in \cref{remark:simplicial-sheaves-terminology} of this current paper.} are model categories, and thus partial model categories.
        \end{definition}

        \begin{definition}
            We begin by formally defining some relative categories (all of which are actually \emph{homotopical} categories).
            All complexes are bounded \emph{cochain} complexes; if we don't say what the morphisms are, then they are simply morphisms of cochain complexes (i.e. degree-wise morphisms that commute with the differentials); and if we don't say what the weak equivalences are, then they are simply quasi-isomorphisms of complexes.
            \begin{itemize}
                \item $\gcohX$ is the category of complexes $C^\bullet$ of sheaves of $\OO_X$-modules such that $C^\bullet$ is quasi-isomorphic to a complex of coherent sheaves.
                \item $\gcohUX$ is the category of complexes $C^\bullet$ of sheaves of $\OO_X$-modules such that the restriction of $C^\bullet$ to any $U\in\cover$ is quasi-isomorphic to a complex of coherent sheaves on $U$, i.e. such that $C^\bullet|U\in\mathsf{Coh}(U)$.
                \item $\gccohX$ is the category of complexes $C^\bullet$ of sheaves of $\OO_X$-modules such that $C^\bullet$ has coherent (internal) cohomology, i.e. such that each $\Ker d_{i}/\Im d_{i-1}$ is a coherent sheaf.
                \item $\cartshX$ is the category of \emph{cartesian} complexes of sheaves on the nerve.
                    Note that morphisms between two such complexes are maps such that, in every \emph{internal} degree (i.e. in every degree \emph{of the complex}), we have a \emph{morphism} of sheaves on the nerve (and such that these commute with the differentials of the complexes); weak equivalences are given by morphisms of complexes such that, in each \emph{simplicial} degree, we have a quasi-isomorphism of complexes.
                \item $\cartvectX$ is the full subcategory of $\cartshX$ consisting of complexes that are \emph{locally (with respect to $\cover$) quasi-isomorphic to} a (cartesian) complex of \emph{locally free} sheaves on the nerve.
                    That is, for $\mathcal{F}^{\bullet,\anotherbullet}\in\cartshX$ to be in $\cartvectX$, there must exist, for all $U_{\alpha_0\ldots\alpha_p}\in\cover$, a complex $\mathcal{G}_{U_{\alpha_0\ldots\alpha_p}}^{\anotherbullet}$ of locally free sheaves on $U_{\alpha_0\ldots\alpha_p}$ such that we have a quasi-isomorphism
                    \[
                        \mathcal{F}^{p,\anotherbullet}_{\alpha_0\ldots\alpha_p} \simeq \mathcal{G}_{U_{\alpha_0\ldots\alpha_p}}^\anotherbullet.
                    \]
                    Similarly, $\cartcohX$ is the full subcategory of $\cartshX$ consisting of complexes that are \emph{locally (with respect to $\cover$) quasi-isomorphic to} a (cartesian) complex of coherent sheaves on the nerve.
                \item $\greenX$ is the full subcategory of $\cartvectX$ spanned by objects that are locally (with respect to $\cover$) quasi-isomorphic to some Green complex (\cref{definition:Gre{}en})
                    The fact that this actually is a subcategory is justified by \cref{corollary:green-gives-cartesian-things}.
                    By definition, every object of $\greenX$ is a Gre{}en complex.
            \end{itemize}

            Note that all of these categories that depend on $\cover$ are natural in the choice of cover: taking a refinement $\anothercover\supset\cover$ induces a functor e.g. $\greenX\to\greenVX$.
            This lets us take homotopy colimits (in the sense of \cref{lemma:filtered-poset-colimit}) over refinements of covers, e.g. $\hocolim_\cover\greenX$.
        \end{definition}

        \begin{remark}\label{remark:subcats-are-locally-q-iso-to}
            Our categories are defined as having ``objects \emph{locally quasi-isomorphic to} a certain class of objects'', and \emph{not} ``objects \emph{in} a certain class of objects''.
            For example, an object $\mathscr{F}^\bullet\in\gcohUX$ is such that, for all $U\in\cover$, there exists some complex of coherent sheaves $\mathscr{G}_U^\bullet$ such that $\mathscr{G}_U^\bullet\simeq \mathscr{F}^\bullet\restricted U$.
            This does \emph{not} a priori imply the existence of some global complex of coherent sheaves $\mathscr{G}^\bullet$ such that $\mathscr{G}^\bullet\simeq\mathscr{F}^\bullet$.
            This subtlety makes the proofs in this section more technical than they morally are.
            %
            % The reason for such subtlety is that we wish to end up with a statement about $\gccohX$, which is already closed under quasi-isomorphisms: if $\mathscr{F}^\bullet$ is in $\gccohX$ and $\mathscr{G}^\bullet$ is quasi-isomorphic to $\mathscr{F}^\bullet$, then $\mathscr{G}^\bullet$ is also in $\gccohX$, since its cohomology is \emph{isomorphic} to that of $\mathscr{F}^\bullet$, which is coherent, by definition.
            % If we wish to have any hope of finding some equivalent category of sufficiently well-behaved sheaves on the nerve then we must expect it to also be closed under quasi-isomorphisms.
            % Aphoristically, if we have a (real life) chain that we wish to make stronger, and we start by strengthening one link, then we must strengthen all other links by the same amount.
        \end{remark}

        \begin{remark}
            Given a full embedding $(\mathcal{C},\mathcal{W}')\hookrightarrow(\mathcal{D},\mathcal{W})$ of one homotopical category into another such that $\mathcal{W}'=\mathcal{C}\cap\mathcal{W}$, there are \emph{two} $(\infty,1)$-categories defined by $\mathcal{C}$ that, in general, do \emph{not} agree:
            \begin{enumerate}
                \item the localisation $\LL_{\mathcal{W}'}\mathcal{C}$ of $\mathcal{C}$ along $\mathcal{W}'$;
                % \item $\LL_{\mathcal{W}'}{\mathcal{C}}$, given by localising $\mathcal{C}$ along $\mathcal{W}'$;
                \item the full sub-$(\infty,1)$-category of $\LL_{\mathcal{W}}\mathcal{D}$ consisting of the objects of $\mathcal{C}$.
                % \item $\LL_\mathcal{D}{\mathcal{C}}$, given by taking the full sub-$(\infty,1)$-category of $\mathcal{D}$ spanned by $\mathcal{C}$.
            \end{enumerate}
            As mentioned in \cref{remark:subcats-are-locally-q-iso-to}, all the subcategories we work with here are defined in the same way: given some $(\mathcal{D},\mathcal{W})$, we construct some $(\mathcal{C},\mathcal{W}')$ by taking the full subcategory of $\mathcal{D}$ of objects that are (locally) connected via $\mathcal{W}$ to objects satisfying some specific property.
            Because of this, we are interested in the \emph{latter} construction (i.e. the full subcategory of $\LL_\mathcal{W}\mathcal{D}$ spanned by $\mathcal{C}$), which we denote by just $\LL\mathcal{C}$, avoiding more suggestive (but lengthier) notation, such as $\LL_{\mathcal{W}',\mathcal{D}}\mathcal{C}$.
        \end{remark}

        \begin{remark}
            In the algebraic setting, there is an equivalence between the category of complexes of coherent sheaves and the category of complexes of sheaves with coherent (internal) cohomology.
            In the analytic case, things are more subtle, and so we have to take a slightly longer route to prove our desired result.
            Indeed, to the best of our knowledge, the question of whether or not $\gcohX$ and $\gccohX$ are equivalent is still open, except in low dimensions, where it is known to be true (see \cite[§2.2.2]{Yu2013}).

            Note that, given a refinement of our cover $\anothercover\supset\cover$, we have full embeddings
            \[
                \gcohX \hookrightarrow \gcohUX \hookrightarrow \gcohVX \hookrightarrow \gccohX
            \]
            that preserve and reflect quasi-isomorphisms.
        \end{remark}

        In summary, we have the following diagram \emph{of $1$-categories}:
        \begin{equation}
        \label{equation:the-big-picture-diagram}
            \begin{tikzcd}[row sep=huge]
                \greenX
                    \ar[r,tail,"\numberincircle{1}"]
                & \cartvectX
                    \ar[r,tail,"\numberincircle{2}"]
                & \cartcohX
                    \ar[r,dashed,"\numberincircle{3}"]
                & \gcohUX
            \end{tikzcd}
        \end{equation}
        where we write $\rightarrowtail$ to denote a fully-faithful functor, and where we will define $\numberincircle{3}$ later, and, even then, only at the level of $(\infty,1)$-categories.
        Note that $\numberincircle{1}$ and $\numberincircle{2}$ are fully-faithful by definition (as full subcategories).

        \bigskip

        \begin{center}
            \emph{Our goal for the rest of this section is to prove that, when we localise all the categories in \cref{equation:the-big-picture-diagram}, all the functors become equivalences of $(\infty,1)$-categories.}

            \emph{This is the content of \cref{theorem:the-main-coherent-theorem,corollary:the-main-coherent-corollary}.}
        \end{center}

        \medskip

        \begin{remark}
            We omit the localisation notation from our functors: we will write $\numberincircle{i}$ to mean both the functor between homotopical categories and the induced functor between their localisations.
        \end{remark}
        
        Our first goal is to show that \cref{equation:the-big-picture-diagram} descends to a diagram at the level of localisations.
        For this we need to know that all our functors $\numberincircle{i}$ really are functors of relative categories, in that they preserve weak equivalences.
        This is automatically true for $\numberincircle{1}$ and $\numberincircle{2}$, since they are inclusions of full subcategories.
        So now we need to construct the functor $\numberincircle{3}$, and we do this by building an adjunction of $(\infty,1)$-categories, step by step.

        Writing $i\colon\nerve{\bullet}\to X$ to mean the map given by inclusion of each open subset into $X$, we have an adjunction\footnote{Note that, to eventually agree with the orientation of our diagram, we write our adjunctions as $(R\vdash L)$ instead of as $(L\dashv R)$.}
        \begin{equation*}
            \shX : (i_* \vdash i^*) : \gshX.
        \end{equation*}
        Then, recalling that the limit functor can be defined as being the right adjoint to the constant diagram functor, and writing $\shX$ to mean the category of complexes of sheaves of $\OO_{\nerve{\bullet}}$-modules, thought of as a category of $\nerve{\bullet}$-diagrams, we can compose this adjunction with the above to get a Quillen adjunction
        \begin{equation*}
            \shX : (\mathrm{lim}\circ i_* \vdash \cst\circ i^*) : \gshX
        \end{equation*}
        where being Quillen follows from the fact that the pullback/pushforward adjunction and the limit/constant adjunction are both Quillen.

        So we are now in the following situation: we have a diagram
        \begin{equation*}
            \begin{tikzcd}
                \shX
                    \ar[r,shift left=2,"\mathrm{lim}\circ i_*"{name=U}]
                & \gshX
                    \ar[l,shift left=2,"\cst\circ i^*"{name=D}]\\
                \cartcohX
                    \ar[u,hook]
                & \gcohUX
                    \ar[u,hook]
                \ar[from=U,to=D,phantom,"\top"]
            \end{tikzcd}
        \end{equation*}
        where the adjunction is Quillen.
        Deriving the functors (and localising the categories) then gives us the diagram
        \begin{equation*}
            \begin{tikzcd}
                \LL{\shX}
                    \ar[r,shift left=2,"\mathbb{R}(\mathrm{lim}\circ i_*)"{name=U}]
                & \LL{\gshX}
                    \ar[l,shift left=2,"\mathbb{L}(\cst\circ i^*)"{name=D}]\\
                \LL{\cartcohX}
                    \ar[u,hook]
                & \LL{\gcohUX}
                    \ar[u,hook]
                \ar[from=U,to=D,phantom,"\top"]
            \end{tikzcd}
        \end{equation*}
        and we wish to know if the adjunction restricts to give an adjunction of the subcategories.
        Following \cite[Lemma~2.2.2.13]{Toen&Vezzosi2008}, we write $\smallint$ to mean the total right derived functor $\mathbb{R}(\mathrm{lim}\circ i_*)$; we write $\iota^*$ to mean the total left derived functor $\mathbb{L}(\cst\circ i^*)$.
        So we want to show that
        \begin{align*}
            \smallint \colon \LL{\cartcohX} &\to \LL{\gcohUX}\\
            \LL{\cartcohX} &\leftarrow \LL{\gcohUX} : \iota^*
        \end{align*}
        (where we omit the restriction of the functors from our notation).

        \begin{lemma}
            The image of $\iota^*\colon\LL{\gcohUX}\to\LL{\shX}$ is contained in $\LL{\cartcohX}$.

            \begin{proof}
                Let $\mathscr{F}^\anotherbullet\in\LL{\gcohUX}$.
                The pullback functor $i^*$ is exact since it is given by the topological pullback tensored with the structure sheaf, and the topological pullback is exact, and tensoring along a disjoint union of open immersions is also exact.
                Then, since the constant diagram functor is (trivially) also exact, we see that $\iota^*$ is just the pullback to the nerve (as with global vector bundles in \cref{example:global-vector-bundles}).
                But, as mentioned there, the simplicial maps are then simply identity maps, which means that the resulting object is indeed cartesian; and being coherent is a local property, so it suffices to check it on each $U_{\alpha_0\ldots\alpha_p}$ in $\nerve{p}$, but, over such an open set, $\iota^*\mathscr{F}^j$ is simply \mbox{$\mathscr{F}^j\restricted U_{\alpha_0\ldots\alpha_p}$}, which is coherent by definition.
            \end{proof}
        \end{lemma}

        \begin{remark}\label{remark:what-is-the-pushforward-as-cosimplicial-object}
            It is a good idea to fully understand the cosimplicial structure of $\iota_*\mathcal{F}^{\bullet,\anotherbullet}$ before proceeding, so we spell out all the details here.

            Recall that $\iota_* = \mathbb{R}i_* = i_*$, where $i\colon\nerve{\bullet}\to X$.
            We want to describe what $(\iota_*\mathcal{F}^{\bullet,\anotherbullet})^p$ is for each $p\in\mathbb{N}$, as well as how these `fit together' to give a cosimplicial object.
            This can be explained by just improving our notation: write $i_p\colon\nerve{p}\to X$ to mean the map given by inclusion of each $U_{\alpha_0\ldots\alpha_p}$ in $\nerve{p}$ into $X$, so that $i$ is exactly the data of $(i_p)_{p\in\mathbb{N}}$.
            Then define
            \[
                (\iota_*\mathcal{F}^{\bullet,\anotherbullet})^p
                =
                (i_p)_*\mathcal{F}^{p,\anotherbullet}
                =
                \bigoplus_{(\alpha_0\ldots\alpha_p)} (i_p)_* \mathcal{F}^{p,\anotherbullet}_{\alpha_0\ldots\alpha_p}.
            \]

            Then, given some $\varphi\colon[p]\to[q]$ in $\Delta$, we want to know how to define
            \[
                (i_p)_*\mathcal{F}^{p,\anotherbullet}
                \xrightarrow{(\iota_*\mathcal{F}^{\bullet,\anotherbullet})^\yetanotherbullet(\varphi)}
                (i_q)_*\mathcal{F}^{q,\anotherbullet},
            \]
            but, using the pullback/pushforward adjunction, this is the same as asking for a map
            \[
                (i_q)^*(i_p)_*\mathcal{F}^{p,\anotherbullet}
                \longrightarrow
                \mathcal{F}^{q,\anotherbullet}.
            \]

            Firstly, we claim, for $p<q$ (dealing with the other case shortly), that there is a natural map $(i_q)^*(i_p)_*\mathcal{F}^{p,\anotherbullet} \rightarrow (\nerve{\bullet}\varphi)^*\mathcal{F}^{p,\anotherbullet}$; secondly, we claim that this gives us the map that we want.
            To see the first claim, we appeal to the geometric nature of pushforwards and pullbacks: since the $U_\alpha$ do not necessarily have trivial intersection with one another, it is not necessarily true that $(i^*i_*\mathcal{F})\restricted U = \mathcal{F}\restricted U$; what \emph{is} true, however, is that the right-hand side is a direct summand of the left:
            \[
                (i^*i_*\mathcal{F})\restricted U
                \cong
                \mathcal{F}\restricted U
                \oplus
                \mathcal{G}.
            \]
            This gives us the first claim (for $p<q$; we deal with the other case shortly), since restriction is the same as pulling back along $\nerve{\bullet}\varphi$.
            For the second claim, by definition of what it means to be an element of $\cartcohX$, we have maps $(\nerve{\bullet}\varphi)^*\mathcal{F}^{p,\anotherbullet}\to\mathcal{F}^{q,\anotherbullet}$ for every $\varphi\colon[p]\to[q]$.
            Combining all of the above then gives us the desired maps, and thus the cosimplicial structure.

            Finally then, when $p>q$ we can do something similar.
            Here the map $\nerve{\bullet}\varphi$ is given by inserting degenerate intersections (obtaining things that look like $U_{\alpha_0\ldots\alpha_i\alpha_i\ldots\alpha_p}$), and so it again isn't necessarily the case that $(i_q)^*(i_p)_*\mathcal{F}^{p,\anotherbullet} = (\nerve{\bullet}\varphi)^*\mathcal{F}^{p,\anotherbullet}$, since we have intersections $U_{\alpha_0\ldots\alpha_p}$ that can be strictly smaller than any $U_{\beta_0\ldots\beta_q}$.
            But we can still construct some $(i_q)^*(i_p)_*\mathcal{F}^{p,\anotherbullet} \to \mathcal{F}^{q,\anotherbullet}$ by precomposing the $(\nerve{\bullet}\varphi)^*\mathcal{F}^{p,\anotherbullet}\to\mathcal{F}^{q,\anotherbullet}$ with the projection maps
            \[
                (i_q)^*(i_p)_*\mathcal{F}^{p,\anotherbullet}
                =
                (i_p)_*\mathcal{F}^{p,\anotherbullet}\restricted U_{\beta_0\ldots\beta_q}
                =
                \bigoplus_{(\alpha_0\ldots\alpha_p)} (i_p)_*\mathcal{F}_{\alpha_0\ldots\alpha_p}^{p,\anotherbullet} \restricted U_{\beta_0\ldots\beta_q}
                \twoheadrightarrow
                (i_p)^*\mathcal{F}_{\beta_0\ldots\beta_q\beta_q\ldots\beta_q}^{q,\anotherbullet}
            \]
            where the first `equality' really means that we work locally over each $U_{\beta_0\ldots\beta_p}$, and we write $(\beta_0\ldots\beta_q\beta_q\ldots\beta_q)$ to mean some degenerate embedding of $(\beta_0\ldots\beta_q)$ into $\nerve{p}$.
        \end{remark}

        \begin{lemma}\label{lemma:image-of-int-is-in-cohU}
            The image of $\smallint\colon\LL{\cartcohX}\to\LL{\gshX}$ is contained in $\LL{\gcohUX}$.

            \begin{proof}
                Before giving the proof, we give a short summary of how it will proceed, to save anybody familiar with such arguments the arduous task of following the notation.
                In particular, this proof is incredibly similar to that of \cref{lemma:counit-is-iso}, which is really an analytic version of \cite[Lemma~2.2.2.13]{Toen&Vezzosi2008} in that it follows the same line of argument.
                We can argue everything locally on some $U\in\cover$; by definition, weak equivalences are quasi-isomorphisms; we can use the total complex construction to calculate the homotopy limit in the definition of $\smallint$; the fact that our complexes are cartesian gives us a weak equivalence between this total complex and the (total complex of the) Čech complex of the simplicial-degree-zero part of our original complex; the latter is weakly equivalent to the simplicial-degree-zero part of our original complex (since all covers can be taken to be Stein); a commuting triangle then tells us that the desired quasi-isomorphism is indeed a quasi-isomorphism.

                \medskip

                Let $\mathcal{F}^{\bullet,\anotherbullet}\in\LL{\cartcohX}$ and $U\in\cover$.
                Then $\mathcal{F}^{0,\anotherbullet}\restricted U$ is a complex of coherent sheaves on $U$, and so it would suffice to show that there is a quasi-isomorphism
                \[
                    \mathcal{F}^{0,\anotherbullet} \restricted U
                    \congto
                    \left(
                        \smallint \mathcal{F}^{\bullet,\anotherbullet}
                    \right) \restricted U
                \]
                in $\gcohUX$.

                Firstly, note that there exists a good candidate morphism: the cosimplicial structure of $\mathcal{F}$ means that we have a morphism
                \[
                    i_p^*\mathcal{F}^{0,\anotherbullet}
                    \to
                    \mathcal{F}^{p,\anotherbullet}
                \]
                of sheaves over $\nerve{p}$ for all $p\in\mathbb{N}$, where $i_p\colon\nerve{p}\to X$;
                by the pull/push adjunction, this gives us a morphism
                \[
                    \mathcal{F}^{0,\anotherbullet}
                    \to
                    (i_p)_*\mathcal{F}^{p,\anotherbullet}
                \]
                and so, by the universal property of the homotopy limit (since $\mathbb{R}\lim=\holim$), we get a zig-zag\footnote{The fact that the zig-zag will consist of three such arrows is mentioned in e.g. \cite[§1.1~(i)]{Barwick&Kan2013}.} of morphisms
                \[
                    \mathcal{F}^{0,\anotherbullet}
                    \congfrom
                    \bullet
                    \to
                    \bullet
                    \congfrom
                    \holim_p(i_p)_*\mathcal{F}^{p,\anotherbullet}
                    =
                    \smallint \mathcal{F}^{\bullet,\anotherbullet}
                \]
                which, writing $\leftrightsquigarrow$ for the zig-zag $\congfrom\bullet\to\bullet\congfrom$, and denoting restriction to $U$ by a subscript $U$, induces
                \begin{equation}\label{equation:candidate-map-to-be-q-iso}
                    \mathcal{F}_U^{0,\anotherbullet}
                    \leftrightsquigarrow
                    \left(
                        \smallint \mathcal{F}^{\bullet,\anotherbullet}
                    \right)_U.
                \end{equation}

                Now we wish to show that \cref{equation:candidate-map-to-be-q-iso} is indeed a quasi-isomorphism.
                We can (justified by \cref{remark:reedy-totalisation-etc}) calculate the homotopy limit with the \define{total construction}: writing $F^{\bullet,\anotherbullet}$ to mean the cosimplicial object $((i_\bullet)_*)\mathcal{F}^{\bullet,\anotherbullet})$, we define
                \begin{align*}
                    \Tot(F)^n &= \bigoplus_{\ell\in\mathbb{N}} F^{\ell,n-\ell}
                \\  \d_{\Tot(F)} &= \d_F + (-1)^n\check\delta
                \end{align*}
                where $\d_F$ is the differential $\d_{\mathcal{F}^{\bullet,\anotherbullet}}$ coming from the $\anotherbullet$-grading of $\mathcal{F}^{\bullet,\anotherbullet}$, and where
                \[
                    \check\delta^n = \sum_{i=0}^{n+1} (-1)^i F^{\bullet,\anotherbullet}(f_{n+1}^i)
                \]
                is the alternating sum of coface maps, whose action is given by the cosimplicial structure of $F^{\bullet,\anotherbullet}$.
                To see that this makes sense in terms of degrees, note that
                \begin{align*}
                    \check{\delta}^m\colon
                    F^{\ell,m}
                    &\to
                    F^{\ell+1,m}
                \\  \d_F^\ell\colon
                    F^{\ell,m}
                    &\to
                    F^{\ell,m+1}.
                \end{align*}
                So \cref{equation:candidate-map-to-be-q-iso} becomes
                \begin{equation}
                    \mathcal{F}_U^{0,\anotherbullet}
                    \leftrightsquigarrow
                    \Tot(F)_U^\anotherbullet,
                \end{equation}
                and to show that this is a quasi-isomorphism, it is enough to show that
                \begin{equation}
                \label{equation:candidate-map-to-be-q-iso-sections}
                    \mathcal{F}_U^{0,\anotherbullet}(V)
                    \leftrightsquigarrow
                    \Tot(F)_U^\anotherbullet(V)
                \end{equation}
                is a quasi-isomorphism for all $V\in\anothercover$, where $\anothercover$ is some cover of $U$.

                The cartesian condition gives us quasi-isomorphisms
                \[
                    \mathcal{F}_{\beta_0\ldots\beta_p}^{0,\anotherbullet}
                    \congto
                    \mathcal{F}_{\beta_0\ldots\beta_p}^{p,\anotherbullet}
                \]
                of complexes of (coherent) sheaves over any $V_{\beta_0\ldots\beta_p}\in\anothercover$, and so we can refine $\anothercover$ to a cover (which we can always take to be Stein) $\yetanothercover=\{W_\gamma\}$ of $V_{\beta_0\ldots\beta_p}$ such that we have quasi-isomorphisms
                \begin{equation}
                \label{equation:q-iso-of-sections}
                    \mathcal{F}_{\beta_0\ldots\beta_p}^{0,\anotherbullet}(W_{\gamma_0\ldots\gamma_p})
                    \congto
                    \mathcal{F}_{\beta_0\ldots\beta_p}^{p,\anotherbullet}(W_{\gamma_0\ldots\gamma_p})
                \end{equation}
                of complexes of abelian groups.

                But, as complexes of sheaves over $V_{\beta_0\ldots\beta_p}$, we trivially have that
                \[
                    \mathcal{F}_{\beta_0\ldots\beta_p}^{p,\anotherbullet}
                    =
                    (i_p)_* \mathcal{F}_{\beta_0\ldots\beta_p}^{p,\anotherbullet} \restricted W_{\beta_0\ldots\beta_p}
                \]
                and so the right-hand side of \cref{equation:q-iso-of-sections} is exactly $F_U^{\bullet,\anotherbullet}(V_{\beta_0\ldots\beta_p})$; further, the left-hand side is exactly $\cech_{\yetanothercover}^\bullet(\mathcal{F}_U^{0,\anotherbullet})$.
                Together, then, this tells us that \cref{equation:q-iso-of-sections} gives a morphism of bicomplexes
                \[
                    \cech_{\yetanothercover}^\bullet(\mathcal{F}_U^{0,\anotherbullet})
                    \to
                    F_U^{\bullet,\anotherbullet}(V_{\beta_0\ldots\beta_p})
                \]
                which is a quasi-isomorphism on each row.
                By a classical spectral sequence argument\footnote{Taking the mapping cone, applying the spectral sequences associated to a bicomplex, and using induction on the number of non-zero rows, combined with the fact that the direct limit functor for complexes of abelian groups is exact.} we can show that such a morphism of bicomplexes induces a quasi-isomorphism on the respective total complexes:
                \[
                    \Tot\cech_{\yetanothercover}^\bullet(\mathcal{F}_U^{0,\anotherbullet})
                    \congto
                    \Tot F_U^{\bullet,\anotherbullet}(V_{\beta_0\ldots\beta_p}).
                \]
                We can see that the triangle
                \[
                    \begin{tikzcd}
                        \mathcal{F}_U^{0,\anotherbullet}(V_{\beta_0\ldots\beta_p})
                            \ar[r,leftrightsquigarrow]
                        & \Tot F_U^{\bullet,\anotherbullet}(V_{\beta_0\ldots\beta_p})
                    \\  \Tot\cech_{\yetanothercover}^\bullet(\mathcal{F}_U^{0,\anotherbullet})
                            \ar[u]
                            \ar[ur]
                    \end{tikzcd}
                \]
                commutes, where the horizontal arrow is exactly the zig-zag \cref{equation:candidate-map-to-be-q-iso-sections}.
                But the vertical arrow is a quasi isomorphism (because taking the Čech complex with respect to a Stein cover gives a resolution), and we have just shown that the diagonal arrow is a quasi-isomorphism (by the cartesian condition), and so the middle arrow in the horizontal zig-zag must also be a quasi-isomorphism.
            \end{proof}
        \end{lemma}

        \medskip

        In summary then, we have an adjunction
        \begin{equation}
        \label{equation:adjunction-of-localisations}
            \LL{\cartcohX}:(\smallint\vdash\iota^*):\LL{\gcohUX}
        \end{equation}

        \begin{definition}\label{definition:defining-numberincircle3}
            $\numberincircle{3}=\smallint$.
        \end{definition}

        \begin{remark}\label{remark:reedy-totalisation-etc}\footnote{We thank Maximilien Péroux for helpful discussions concerning this remark.}
            The projective model structure on non-negatively-graded cochain complexes gives us a simplicial model category (by the dual of Dold-Kan), and so, if we can show that $F^{\yetanotherbullet,\anotherbullet}$ is Reedy fibrant, then we can apply \cite[Theorem~19.8.7]{Hirschhorn2003}, which tells us that $\holim F^{\yetanotherbullet,\anotherbullet}\simeq\Tot(F^{\yetanotherbullet,\anotherbullet})$, for some abstract definition of $\Tot$.
            The fact that the totalization (in the sense of Hirschhorn) agrees with the totalization of a bicomplex (in the usual homological algebra sense), and that the Bousfield-Kan spectral sequence and the spectral sequence(s) associated to a bicomplex coincide, can be found in \cite[{}III.1.1.13]{Fresse2017}.

            To show Reedy fibrancy of some cosimplicial object $X^\bullet$, we need to show that the maps $X^n\to M_n(X^\bullet)$ are fibrant for all $n\geqslant0$, where $M_n(X^\bullet)$ is the {matching object} given by
            \[
                M_n(X^\bullet)
                =
                \lim_{\substack{\varphi\colon[n]\twoheadrightarrow[i] \\ i\neq n}} X^i.
            \]
            We can write this down more explicitly as
            \begin{align*}
                M_0(X^\bullet)
                &=
                \{*\}
            \\  M_1(X^\bullet)
                &=
                X^0
            \\  M_2(X^\bullet)
                &=
                X^1 \times_{X^0} X^1
            \end{align*}
            and so on.

            It is a purely formal consequence of the simplicial identities that any \emph{simplicial} object in a model category whose cofibrations are the monomorphisms (such as the category of complexes with the injective model structure) is Reedy \emph{co}fibrant;
            formally dual to this is the fact that any \emph{co}simplicial object in a model category whose fibrations are the epimorphisms (such as the category of complexes with the projective model structure) is Reedy \emph{fibrant}.

            For example, the fact that $X^1\to M_1(X^\bullet)=X^0$ is fibrant (i.e. an epimorphism) is due to the fact that it admits a right inverse (namely either of the two face maps $X^0\to X^1$), thanks to the (co)simplicial identities (namely $s_0^0 f_1^0 = s_0^0 f_0^0 = \id_{[0]}$).
        \end{remark}

    \subsection{Equivalences}

        \begin{lemma}\label{lemma:left-adjoint-is-conservative}
            The left adjoint $\iota^*$ of the adjunction~(\ref{equation:adjunction-of-localisations}) is conservative.
            \begin{proof}
                Let $f\colon \mathscr{F}^\anotherbullet\to\mathscr{G}^\anotherbullet$ be a morphism in $\LL{\gccohX}$ such that $\iota^*f$ is an isomorphism in $\LL{\cartcohX}$.
                By definition of the weak equivalences in $\cartcohX$, and the calculation of $\iota^*$ in \cref{definition:defining-numberincircle3}, over each $U_{\alpha_0\ldots\alpha_p}\in\nerve{p}$, this says that the map
                \[
                    \iota^*f = f\restricted U_{\alpha_0\ldots\alpha_p}
                    \colon
                    \mathscr{F}^\anotherbullet \restricted U_{\alpha_0\ldots\alpha_p}
                    \to
                    \mathscr{G}^\anotherbullet \restricted U_{\alpha_0\ldots\alpha_p}
                \]
                is a quasi-isomorphism.
                But this is just saying that every restriction of $f$ to some open subset of $X$ is a quasi-isomorphism, which implies that $f$ is a quasi-isomorphism, i.e. a weak equivalence in $\gccohX$.
            \end{proof}
        \end{lemma}

        \begin{lemma}\label{lemma:weak-equivalences-in-cartshX}
            Let $\varphi^\bullet\colon\mathcal{F}^\bullet\to\mathcal{G}^\bullet$ be a morphism\footnote{We omit from our notation the internal grading $\anotherbullet$ of the complexes, writing e.g. $\mathcal{F}^\bullet$ instead of $\mathcal{F}^{\bullet,\anotherbullet}$.} in $\cartshX$.
            Then $\varphi^\bullet$ is a weak equivalence if and only if $\varphi^0\colon\mathcal{F}^0\to\mathcal{G}^0$ is a weak equivalence.
            \begin{proof}
                If $\varphi^\bullet$ is a weak equivalence then, by definition, each $\varphi^p$ is a weak equivalence, and so it remains only to show `if' part of the claim.
                Recalling that pulling back along the inclusion of an open subset is the same as restriction to that same open subset, the cartesian condition on $\mathcal{F}^\bullet$ tells us that, for all $U_{\alpha_0\ldots\alpha_p}\in\cover$,
                \[
                    \mathcal{F}^0 \restricted U_{\alpha_0\ldots\alpha_p}
                    \congto
                    \mathcal{F}^p \restricted U_{\alpha_0\ldots\alpha_p}
                \]
                as (complexes of) sheaves over $U_{\alpha_0\ldots\alpha_p}$, and similarly for $\mathcal{G}^\bullet$.
                Combining this with the commutative square that follows from the definition of what it means to be a morphism of sheaves on the nerve, we get the commutative square
                \[
                    \begin{tikzcd}
                        \mathcal{F}^0
                            \ar[r,"\sim"]
                            \ar[d,swap,"\varphi^0"]
                       &\mathcal{F}^p
                            \ar[d,"\varphi^p"]
                    \\  \mathcal{G}^0
                            \ar[r,"\sim"]
                       &\mathcal{G}^p
                    \end{tikzcd}
                \]
                from which it follows that, if $\varphi^0\colon\mathcal{F}^0\to\mathcal{G}^0$ is a weak equivalence, then so too is $\varphi^p$ for all $p\in\mathbb{N}$, and hence also $\varphi^\bullet$.
            \end{proof}
        \end{lemma}

        \begin{lemma}\label{lemma:counit-is-iso}
            The counit
            \[
                \iota^*\circ\smallint\implies\id_{\LL{\cartcohX}}
            \]
            of the adjunction~(\ref{equation:adjunction-of-localisations}) is a weak equivalence.
            \begin{proof}
                Let $\mathcal{F}^{\bullet,\anotherbullet}\in\cartcohX$.
                First of all, by \cref{lemma:weak-equivalences-in-cartshX}, we know that it suffices to show that the counit is a weak equivalence in simplicial degree zero.
                But we can further simplify things: recalling the definition of $\smallint$, and using the fact that $\iota^*$ is simply the pullback along $i$, it suffices to show that the induced morphism
                \begin{equation}
                \label{equation:map-for-counit-lemma}
                    \left(
                        i^* \holim_{p\in\mathbb{N}}(\iota_*\mathcal{F}^{\bullet,\anotherbullet})^p
                    \right) \restricted U_\alpha
                    \longrightarrow
                    \mathcal{F}^{0,\anotherbullet} \restricted U_\alpha
                \end{equation}
                in $\cartcohX$ is a weak equivalence for all $U_\alpha\in\cover$.
                But since the composite $U_\alpha\hookrightarrow\cover\to X$ is exactly $U_\alpha\hookrightarrow X$, we see that pulling back along $i$ and then restricting to $U_\alpha$ is the same as restricting directly to $U_\alpha$.
                Finally, just to simplify notation, we write $F^{\bullet,\anotherbullet}$ to mean the cosimplicial object $(\iota_*\mathcal{F}^{\bullet,\anotherbullet})^\bullet$, we write $\d_F$ to mean the differential $\d_{\mathcal{F}^{\bullet,\anotherbullet}}$ coming from the $\anotherbullet$-grading of $\mathcal{F}^{\bullet,\anotherbullet}$, and we write a subscript $\alpha$ to denote restriction to $U_\alpha$.
                All together then, \cref{equation:map-for-counit-lemma} becomes
                \begin{equation}
                \label{equation:simplified-map-for-counit-lemma}
                    \left(
                        \holim_{p\in\mathbb{N}}F^{p,\anotherbullet}
                    \right) \restricted U_\alpha
                    \longrightarrow
                    \mathcal{F}_\alpha^{0,\anotherbullet}
                \end{equation}

                But now we can proceed almost exactly as in the proof of \cref{lemma:image-of-int-is-in-cohU}: we use the total construction, and construct the same commuting triangle but with the horizontal arrow going in the other direction.
            \end{proof}
        \end{lemma}

        \begin{remark}\label{remark:counit-iso-iff-right-adjoint-fully-faithful}
            The counit of an adjunction being a weak equivalence (or an isomorphism, in the 1-categorical case) is equivalent to the right adjoint being fully faithful.
        \end{remark}

        \begin{lemma}\label{lemma:adjunction-conservative-fully-faithful-is-equiv}
            Let $\mathcal{C}:(L\dashv R):\mathcal{D}$ be an adjunction with $L$ conservative and $R$ fully faithful.
            Then $(L\dashv R)$ gives an equivalence $\mathcal{C}\simeq\mathcal{D}$.
            \begin{proof}
                It suffices to show that $R$ is essentially surjective, so let $c\in\mathcal{C}$, and define $d=L(c)$.
                Then $LR(d)\to d$ is an equivalence (because $R$ being fully faithful is equivalent to the counit of the adjunction being an equivalence).
                But $LR(d)\to d$ is, by definition, $LRL(c)\to L(c)$, and since this is an equivalence and $L$ is conservative, we see that $RL(c)\to c$ is an equivalence. That is, $R(d)\simeq c$.
            \end{proof}
        \end{lemma}

        \begin{corollary}\label{corollary:adjunction-is-equivalence}
            The adjunction~(\ref{equation:adjunction-of-localisations}) gives an equivalence of $(\infty,1)$-categories
            \[
                \LL{\cartcohX} \simeq \LL{\gcohUX}
            \]
            and thus an equivalence
            \[
                \hocolim_\cover\LL{\cartcohX} \simeq \hocolim_\cover\LL{\gcohUX}.
            \]
        \end{corollary}

        \begin{lemma}\label{lemma:green-gives-essential-surjectivity}
            The composite functor
            \begin{equation*}
                \hocolim_\cover\LL{\greenX}
                \xrightarrow{\numberincircle{3}\numberincircle{2}\numberincircle{1}}
                \hocolim_\cover\LL{\gcohUX}
            \end{equation*}
            is essentially surjective.
            \begin{proof}
                Since $\numberincircle{3}$ is an equivalence (\cref{corollary:adjunction-is-equivalence}), it suffices to show that
                \[
                    \numberincircle{2}\numberincircle{1}\colon
                    \hocolim_\cover\LL{\greenX}
                    \to
                    \hocolim_\cover\LL{\cartcohX}
                \]
                is essentially surjective.
                So let $\mathscr{F}^{\bullet,\anotherbullet}\in\cartcohX$.

                By definition, for all $U_\alpha\in\cover$, there exists some complex $\mathscr{G}_{\alpha}^\anotherbullet$ of coherent sheaves on $U_\alpha$ such that $\mathscr{F}^{0,\anotherbullet}_\alpha \simeq \mathscr{G}_\alpha^\anotherbullet$.
                We know that we can always locally resolve $\mathscr{G}_\alpha^\anotherbullet$ by locally free sheaves, and so, by possibly taking a refinement $\anothercover\supset\cover$ (and using $\alpha$, $\beta$, $\ldots$ to now label the open sets $V_\alpha$, $V_\beta$, $\ldots$ of the refinement), we can obtain some (bounded) complex $\mathscr{H}_\alpha^\anotherbullet$ of free sheaves (of finite rank) on $V_\alpha$ such that $\mathscr{F}_\alpha^{0,\anotherbullet} \simeq \mathscr{H}_\alpha^\anotherbullet$.

                But this is simply saying that $\mathscr{F}_\alpha^{0,\anotherbullet}$ is {perfect}, and so \cite[Proposition~1.2.3]{OBrian&etal1985} (or \cite[Proposition~3.20]{Wei2016}), tells us that, after possibly taking another refinement of our cover, there exists some holomorphic twisting resolution of $\mathscr{F}^{0,\anotherbullet}$.
                Applying the construction of Green's resolution then gives us some $\mathcal{E}^{\bullet,\anotherbullet}\congto j^*(\mathscr{F}^{0,\anotherbullet})^\bullet$, where $\mathcal{E}^{\bullet,\anotherbullet}$ is a complex of locally free sheaves on the Čech nerve of $\nerve{0}=\coprod_\alpha V_\alpha$, and $j$ is the map from the Čech nerve of $\nerve{0}$ to $\nerve{0}$ itself.
                Note, however, that the Čech nerve of $\nerve{0}$ is identical to the Čech nerve of $X$, and so it suffices to prove that $j^*(\mathscr{F}^{0,\anotherbullet})^\bullet\simeq\mathscr{F}^{\bullet,\anotherbullet}$ as sheaves on $\nerve{\bullet}$, since $\numberincircle{1}$ and $\numberincircle{2}$ are both simply inclusions of full subcategories.

                By \cref{lemma:weak-equivalences-in-cartshX}, it suffices to show that we have a weak equivalence in simplicial degree zero, but $j^*(\mathscr{F}^{0,\anotherbullet})^0\congto\mathscr{F}^{0,\anotherbullet}$ is simply the identity map.
            \end{proof}
        \end{lemma}

        \begin{lemma}\label{lemma:filtered-poset-colimit}
            Let $\mathcal{C}$ be a partial model category, and $(\mathcal{D}_\lambda)_{\lambda\in\mathcal{P}}$ be a diagram of full subcategories of $\mathcal{C}$ indexed by some filtered poset $\mathcal{P}$.
            Assume further that each $\mathcal{D}_\lambda$ is stable under weak equivalences.
            Then
            \[
                \hocolim_{\lambda\in\mathcal{P}} \LL{\mathcal{D}_\lambda}
                \simeq
                \bigcup_{\lambda\in\mathcal{P}} \LL{\mathcal{D}_\lambda}
            \]
            where $\bigcup_{\lambda\in\mathcal{P}} \LL{\mathcal{D}_\lambda}$ is the full sub-$(\infty,1)$-category of $\LL{\mathcal{C}}$ spanned by the union of the objects of all of the $\LL{\mathcal{D}_\lambda}$.

            In particular, the induced map
            \[
                \hocolim_{\lambda\in\mathcal{P}} \LL{\mathcal{D}_\lambda}
                \to
                \LL{\mathcal{C}}
            \]
            is fully faithful.

            \begin{proof}
                Since partial model categories present $(\infty,1)$-categories, it suffices to prove the corresponding claim for $\hocolim\mathcal{D}_\lambda$ instead of $\hocolim\LL{\mathcal{D}_\lambda}$.
                % Using the Quillen equivalence between relative categories and complete Segal spaces (\cite[Theorem~6.1]{Barwick&Kan2012}) we know that a morphism of relative categories is a cofibration if it is sent to a cofibration of complete Segal spaces.
                Write $Y_\bullet^\lambda$ to mean the complete Segal space $\LL{\mathcal{D}_\lambda}$, and $X_\bullet$ to mean the complete Segal space $\LL{\mathcal{C}}$.
                But since each $\mathcal{D}_\lambda$ is a full subcategory of $\mathcal{C}$ stable under weak equivalences, each $Y_n^\lambda$ is a union of connected components of $X_n$, corresponding to the span of the objects of $\mathcal{D}_\lambda$.
                This means that, in particular, for each $\mathcal{D}_\lambda\hookrightarrow\mathcal{D}_\mu$, the maps $Y_n^\lambda\hookrightarrow Y_n^\mu$ are all closed embeddings, and thus cofibrations.
                Hence
                \[
                    \hocolim_\lambda Y_\bullet^\lambda \cong \colim_\lambda Y_\bullet^\lambda.
                \]
                Similarly, each $Y_n^\lambda\hookrightarrow X_n$ is a closed embedding, and thus a cofibration.

                We claim that
                \[
                    \colim_\lambda Y_\bullet^\lambda \cong \bigcup_\lambda Y_\bullet^\lambda
                \]
                where $\bigcup_\lambda Y_\bullet^\lambda$ is the subspace of $X_\bullet$ spanned by the connected components of all the $Y_\bullet^\lambda$.
                Note that this is the complete Segal space $\bigcup_\lambda\LL{\mathcal{D}_\lambda}$, and so proving the above claim will finish the proof.

                By the universal property of the colimit, we have a commutative diagram
                \[
                    \begin{tikzcd}
                        (Y_\bullet^\lambda)_{\lambda\in\mathcal{P}}
                            \ar[d,swap,"(f_\lambda)_{\lambda\in\mathcal{P}}"]
                            \ar[dr,"\iota"]
                    \\  \colim_\lambda Y_\bullet^\lambda
                            \ar[r,swap,"f"]
                        & \bigcup_\lambda Y_\bullet^\lambda
                    \end{tikzcd}
                \]
                and we wish to show that $f$ is an isomorphism.
                Since $\iota$ is simply the inclusion of each $Y_\bullet^\lambda$ into the union, it is surjective.
                Thus, given any $y\in\bigcup_\lambda Y_\bullet^\lambda$, there exists some $\lambda$ such that $y\in Y_\bullet^\lambda$, and then, by commutativity, $f(f_\lambda(y))=y$, which shows surjectivity.
                To show injectivity, let $y,z\in\colim_\lambda Y_\bullet^\lambda$ be such that $f(y)=f(z)$.
                Since $\mathcal{P}$ is filtered, there exists some $\lambda$ such that $f(y),f(z)\in Y_\bullet^\lambda$.
                Now $f(y)=ff_\lambda f(y)$ and $f(z)=ff_\lambda f(z)$, whence $ff_\lambda f(y)=ff_\lambda f(z)$.
                But $f$ is surjective, and so $ff_\lambda(y)=ff_\lambda(z)$, but $ff_\lambda=\iota_\lambda$ is simply the inclusion of $Y_\bullet^\lambda$ into the union, whence $y=z$.
            \end{proof}
        \end{lemma}

        \begin{lemma}\label{lemma:green-gives-fully-faithful}
            The composite functor
            \begin{equation*}
                \hocolim_\cover\LL{\greenX}
                \xrightarrow{\numberincircle{3}\numberincircle{2}\numberincircle{1}}
                \hocolim_\cover\LL{\gcohUX}
            \end{equation*}
            is fully faithful.
            \begin{proof}
                By \cref{lemma:counit-is-iso,remark:counit-iso-iff-right-adjoint-fully-faithful}, $\numberincircle{3}$ is fully faithful; since $\numberincircle{1}$ and $\numberincircle{2}$ are inclusions of full subcategories, they remain fully faithful at the level of localisations, and we can use \cref{lemma:filtered-poset-colimit} to see that they remain fully faithful after taking homotopy colimits; this means that the composite functor $\numberincircle{3}\numberincircle{2}\numberincircle{1}$ is also fully faithful.
            \end{proof}
        \end{lemma}

        \begin{corollary}
        \label{corollary:equivalence-without-connections}
            There is an equivalence of $(\infty,1)$-categories
            \begin{equation*}
                \hocolim_\cover\LL{\greenX} \simeq \hocolim_\cover\LL{\gcohUX}.
            \end{equation*}
            \begin{proof}
                \cref{lemma:green-gives-essential-surjectivity} tells us that $\numberincircle{3}\numberincircle{2}\numberincircle{1}$ is essentially surjective; \cref{lemma:green-gives-fully-faithful} tells us that it is fully faithful.
                Thus the composite functor $\numberincircle{3}\numberincircle{2}\numberincircle{1}$ is an equivalence, and, since all of the constitutive functors are fully faithful, each one must itself also be an equivalence.
            \end{proof}
        \end{corollary}

\section{Simplicial connections}\label{section:simplicial-connections}

    \begin{definition}
        We introduce the notation $\pi_{p,q}\colon\nerve{p}\times\Delta^q\to\nerve{p}$ for the projection map; we write $\pi_p$ to mean $\pi_{p,p}$.
        Given some vector bundle $\mathcal{E}^\bullet$ on the nerve, we sometimes write $\overline{\mathcal{E}^p}$ to mean $\pi_p^*\mathcal{E}^p$.
    \end{definition}

    \subsection{Connections and true morphisms}

        \begin{definition}
            The \define{Atiyah exact sequence} (or \define{jet sequence}) of a holomorphic vector bundle $E$ on $X$ is the short exact sequence of $\OO_X$-modules
            \[
                0
                \to
                E\otimes\Omega_X^1
                \to
                J^1(E)
                \to
                E
                \to
                0
            \]
            where $J^1(E)=(E\otimes\Omega_X^1)\oplus E$ as a sheaf of $\mathbb{C}_X$-modules, but with an $\OO_X$-action defined by
            \[
                f(s\otimes\omega,t)
                =
                (fs\otimes\omega + t\otimes\d f, ft).
            \]
        \end{definition}

        \begin{definition}
            A \define{holomorphic (Koszul) connection} $\nabla$ on a holomorphic vector bundle $E$ on $X$ is a (holomorphic) splitting of the Atiyah exact sequence of $E$.
            By enforcing the \define{Leibniz rule}
            \[
                \nabla(s\otimes\omega)
                =
                \nabla s\wedge\omega + s\otimes\d\omega
            \]
            we can extend any connection $\nabla\colon E\to E\otimes\Omega_X^1$ to a map $\nabla\colon E\otimes\Omega_X^r\to E\otimes\Omega_X^{r+1}$.
            (Using the same symbol $\nabla$ to denote the connection as well as any such extension is a common abuse of notation.)
            The \define{curvature} $\kappa(\nabla)$ of a connection $\nabla$ is the $\OO_X$-linear map
            \[
                \kappa(\nabla)=\nabla\circ\nabla\colon
                E\to E\otimes\Omega_X^2.
            \]
        \end{definition}

        \begin{definition}\label{definition:true-morphism}
            A morphism $f\colon(A,\nabla_A)\to(B,\nabla_B)$ of vector bundles on $X$ with connections is said to be a \define{true\footnote{These are sometimes called \define{flat} morphisms, but we opt for `true' to avoid overloading the meaning of the word `flat'.} morphism} if
            \[
                \nabla_B\circ f = (f\otimes\id)\circ\nabla_A.
            \]

            Note that, by the Leibniz rule, if such a morphism $f$ is a true morphism, then the morphism $f\colon(A,\nabla_A^r)\to(B,\nabla_B^r)$ is also `true', in some sense: it satisfies
            \[
                \nabla_B^r \circ f = (f\otimes\id) \circ \nabla_A^r.
            \]
            In particular, if $f$ is a true morphism, then $f\colon(A,\kappa(\nabla_A))\to(B,\kappa(\nabla_B))$ is also a `true' morphism (in this more general sense).
        \end{definition}

        \begin{definition}\label{definition:comparison-map}
            Given a vector bundle on the nerve $\mathcal{E}^\bullet$ we define its \define{$i$-th comparison map}
            \[
                \comparison{p}^i(\mathcal{E}^\bullet)
                \colon
                \left(\nerve{\bullet} f_p^i\times\id\right)^*\overline{\mathcal{E}^{p-1}}
                \longrightarrow
                \left(\id\times f_p^i\right)^*\overline{\mathcal{E}^p},
            \]
            for $i\in\{0,\ldots,p-1\}$, to be the map
            \[
                \pi_{p,p-1}^*\left(\nerve{\bullet} f_p^i\right)^*\mathcal{E}^{p-1}
                \xrightarrow{\pi_{p,p-1}^*\left( \mathcal{E}^\bullet f_p^i \right)}
                \pi_{p,p-1}^*\mathcal{E}^p,
            \]
            where we use the fact that
            \begin{align*}
                \left(\nerve{\bullet} f_p^i\times\id\right)^*\overline{\mathcal{E}^{p-1}}
            &   =
                \pi_{p,p-1}^*\left(\nerve{\bullet} f_p^i\right)^*\mathcal{E}^{p-1}
            \\  \left(\id\times f_p^i\right)^*\overline{\mathcal{E}^p}
            &   =
                \pi_{p,p-1}^*\mathcal{E}^p.
            \end{align*}

            If there is no chance of confusion, we usually omit the dependence on $\mathcal{E}^\bullet$ from the notation.
            We often say `the' comparison map $\comparison{p}$ when we really mean `all' comparison maps $\comparison{p}^i$ for $i\in\{0,\ldots,p-1\}$.
        \end{definition}

        \begin{corollary}\label{corollary:Green-has-injective-comparisons}
            Let $\mathcal{E}^{\bullet,\anotherbullet}$ be a Gre{}en complex.
            Then the comparison maps $\comparison{p}^i(\mathcal{E}^{\bullet,\anotherbullet})$ are injective.

            \begin{proof}
                Pulling back along $\pi_{p,p-1}$ is an exact functor, and being Gre{}en tells us that the $\mathcal{E}^{\bullet,\anotherbullet} f_p^i$ are injective.
            \end{proof}
        \end{corollary}

    \subsection{The motivating example}\label{subsection:the-motivating-example}

        In an effort to motivate the definitions in this chapter, we start with a simplified example of what we wish to study: we replace vector bundles with vector spaces, and we replace curvatures of connections with endomorphisms.

        \begin{definition}
            Let $\mathcal{C}$ be the category whose objects are pairs $(V,\varphi)$ of finite-dimensional vector spaces $V$ and endomorphisms $\varphi$, and whose morphisms $f\colon(V,\varphi)\to(W,\psi)$ are the morphisms $f\colon V\to W$ of vector spaces such that $f\circ\varphi=\psi\circ f$.

            Let $E\colon\mathcal{C}\to\mathcal{C}$ be the endofunctor that sends $(V,\varphi)$ to $(V/\Ker\varphi,\varphi)$.
            We define the wide subcategory $\mathcal{W}$ to consist of all morphisms that become isomorphisms after applying $E$, and denote by $\LL_E\mathcal{C}$ the localisation of $\mathcal{C}$ along $\mathcal{W}$.
        \end{definition}

        Recall that the \define{Grothendieck group} $K(\mathcal{C})$ of $\mathcal{C}$ is the group whose elements are isomorphism classes $[A]$ of objects $A\in\mathcal{C}$, and where, for each short exact sequence $0\to A\to B\to C\to 0$ in $\mathcal{C}$, we introduce the relation $[A]-[B]+[C]=0$ in $K(\mathcal{C})$, whence the group operation is given by
        \[
            [(V,\varphi)] + [(W,\psi)]
            :=
            (V\oplus W,\varphi\oplus\psi).
        \]

        \begin{definition}
            An object $(V,\varphi)\in\mathcal{C}$ is \define{flat} if $\varphi=0$.
            We define an equivalence relation $\sim$ on $K(\mathcal{C})$ by saying that flat objects are equivalent to the zero object, i.e. the equivalence relation is generated by
            \[
                [(V,0)]\sim[(0,0)].
            \]
        \end{definition}

        \begin{definition}\label{definition:admissible-morphism-in-C}
            A morphism $f\colon(V,\varphi)\to(W,\psi)$ in $\mathcal{C}$ is \define{admissible} if there exist sub-bundles $V_1\hookrightarrow V$ and $W_1\hookrightarrow W$ such that
            \begin{enumerate}
                \item $V_1\subseteq\Ker\varphi$ and $W_1\subseteq\Ker\psi$;
                \item $f$ restricts to a morphism $V_1\to W_1$;
                \item $f$ descends to an isomorphism $V/V_1\congto W/W_1$.
            \end{enumerate}
        \end{definition}

        \begin{lemma}\label{lemma:admissible-example-definition}
            A morphism $f\colon(V,\varphi)\to(W,\psi)$ in $\mathcal{C}$ is in $\mathcal{W}$ if and only if it is admissible.
            \begin{proof}
                Let $f\colon(V,\varphi)\to(W,\psi)$ be a morphism in $\mathcal{W}$, so that $E(f)$ is an isomorphism.
                Then we take $V_1:=\Ker\varphi$ and $W_1:=\Ker\psi$.

                Conversely, let $f\colon(V,\varphi)\to(W,\psi)$ be admissible.
                Then $V_1\leqslant\Ker\varphi$, and so, by the third isomorphism theorem, $(V/V_1)/(\Ker\varphi/V_1)\cong V/\Ker\varphi$.
                It remains then to show that $E(f)$ restricts to an isomorphism $\Ker\varphi/V_1\congto\Ker\psi/W_1$.
                But $f\circ\varphi=\psi\circ f$, and so $E(f)\circ E(\varphi)\circ E(f)^{-1}=E(\psi)$, whence
                \[
                    E(f)\colon\Ker\varphi/V_1=\Ker E(\varphi)\congto\Ker E(\psi)=\Ker\psi/W_1.\qedhere
                \]
            \end{proof}
        \end{lemma}

        \begin{lemma}\label{lemma:admissible-example-isomorphism}
            There is a (canonical) isomorphism $K(\LL_E\mathcal{C})\congto K(\mathcal{C})/\!\!\sim$.
            \begin{proof}
                It suffices to take the `identity' map, as follows. Take two isomorphic objects $(V,\varphi)\cong(W,\psi)\in \LL_E\mathcal{C}$. 
                Then the isomorphism between them is given either by an isomorphism in $\mathcal{C}$ or by some morphism in $\mathcal{W}$.
                In the former case, we are done; in the latter case, \cref{lemma:admissible-example-definition} tells us that there is an admissible decomposition $f\colon V\cong V_1\oplus V_2\to W_1\oplus W_2\cong W$.
                But then both $(V_1,\varphi)$ and $(W_1,\psi)$ are equivalent to zero, whence $f\colon(V_2,\varphi)\congto(W_2,\psi)$ gives us an isomorphism in $K(\mathcal{C})/\!\!\sim$.

                Conversely, take two isomorphic objects $(V,\varphi)\cong(W,\psi)\in\mathcal{C}$. Then they are also isomorphic in $\LL_E\mathcal{C}$, since all isomorphisms are in $\mathcal{W}$.
                Further, if \mbox{$[(V,\varphi)]\sim[(0,0)]$} in $K(\mathcal{C})/\!\!\sim$, then $\varphi=0$, whence $V/\Ker\varphi=0$, and so $[(V,\varphi)]=[(0,0)]\in K(\LL_E\mathcal{C})$.
            \end{proof}
        \end{lemma}

        Now assume that we have some `invariant polynomial'\footnote{That is, invariant under a change of basis (the $\mathrm{GL}_n$-action), but with the subtlety that we actually need a \emph{sequence} of such polynomials, indexed by $\mathbb{N}$: one for each possible dimension of the vector space of an element of $\mathcal{C}$. We describe these things in more detail when we need them: in the sequel to this paper.} $P$ from $\mathcal{C}^{\otimes n}$ to some (additive, say) abelian group $G$ (such as $\mathbb{C}$).
        If $P$ is additive then it will descend to a well-defined polynomial on $K(\mathcal{C})^{\otimes n}$.
        If $P$ further sends an $n$-fold tensor product of flat objects to zero, then it also descends to a well-defined polynomial on $K(\mathcal{C})/\!\!\sim$.
        Thus, by \cref{lemma:admissible-example-isomorphism}, $P$ is well defined on $K(\LL_E\mathcal{C})$; thus, by \cref{lemma:admissible-example-definition}, the resulting characteristic classes (that is, the values of $P$) are invariant under admissible morphisms.

    \subsection{Admissibility}

        \begin{definition}
            An \define{endomorphism-valued simplicial $r$-form $\omega_\bullet$} on a vector bundle on the nerve $\mathcal{E}^\bullet$ is a family of forms $\omega_\bullet=\{\omega_p\}_{p\in\mathbb{N}}$, where $\omega_p$ is a global section of the sheaf
            \begin{equation*}
                \sheafend\left(\overline{\mathcal{E}^p}\right)\otimes_{\OO_{\nerve{p}}}\Omega_{\nervesimplex{p}}^r
            \end{equation*}
            such that
            \begin{gather*}
                (\comparison{p}^i\otimes\id)
                \circ
                \big( (\nerve{\bullet} f_p^i\times\id)^*\otimes (\nerve{\bullet} f_p^i\times\id)^* \big)
                \,\omega_{p-1}
            \\  = \big( (\id\times f_p^i)^*\otimes(\id\times f_p^i)^* \big)
                \,\omega_p
                \circ
                (\comparison{p}^i\otimes\id)
            \end{gather*}
            as global sections of $\sheafhom\left((\nerve{\bullet} f_p^i\times\id)^*\overline{\mathcal{E}^{p-1}},\,\,(\id\times f_p^i)^*\overline{\mathcal{E}^p}\right)\otimes\Omega_{\nerve{p}\times\Delta^{p-1}}^r$.
        \end{definition}

        \begin{definition}
            We say that an endomorphism-valued simplicial $r$-form $\omega_\bullet$ is \define{admissible} if it is fibrewise $\comparison{p}$-admissible: for all $p\in\mathbb{N}$, all $x\in X$, all $v_x\in\bigwedge^r T_x X$, the endomorphism
            \[
                \overline{\mathcal{E}^p} \restricted \{x\}
                \xrightarrow{\omega_p(v_x)}
                \overline{\mathcal{E}^p} \restricted \{x\}
            \]
            is such that, for all $i\in\{0,\ldots,p-1\}$, the comparison map $\comparison{p}^i(\mathcal{E}\restricted\{x\})$ is admissible (in the sense of \cref{definition:admissible-morphism-in-C}) with respect to the endomorphisms $\omega_p(v_x)$; in other words, if the induced map
            \[
                \left(\nerve{\bullet} f_p^i\times\id\right)^*
                \bigg(
                    \left(\overline{\mathcal{E}^{p-1}}\restricted \{x\}\right)/\Ker\omega_{p-1}(v_x)
                \bigg)
                \xrightarrow{\comparison{p}^i}
                \left(\id\times f_p^i\right)^*
                \bigg(
                    \left(\overline{\mathcal{E}^p}\restricted \{x\}\right)/\Ker\omega_p(v_x)
                \bigg)
            \]
            is an isomorphism.
        \end{definition}

        \begin{lemma}\label{lemma:admissible-endomorphism-valued-form-criterion}
            For an endomorphism-valued simplicial form $\omega_\bullet$ on $\mathcal{E}^\bullet$ to be admissible it is sufficient to ask, for all $p\in\mathbb{N}$, for sub-bundles $L^p,M^p\hookrightarrow\overline{\mathcal{E}^p}$, lying in the kernel of the endomorphism part of the $\omega_p$, such that the comparison maps $\comparison{p}^i$ restrict to isomorphisms
            \[
                \comparison{p}^i\colon
                \left(\nerve{\bullet} f^i_p\times\id\right)^*\left(\overline{\mathcal{E}^{p-1}}/L^{p-1}\right)
                \congto
                \left(\id\times f_p^i\right)^*\left(\overline{\mathcal{E}^p}/M^p\right).
            \]
            \begin{proof}
                If $L^p$ lies in the kernel of (the endomorphism part of) $\omega_p$ then, in particular, $L^p\restricted \{x\}$ lies in the kernel of $\omega_p(v_x)$ for any $v_x\in\bigwedge^r T_x X$ (and similarly for $M^p$).
                Then we appeal to \cref{lemma:admissible-example-definition}.
            \end{proof}
        \end{lemma}

    \subsection{Simplicial connections}

        \begin{definition}
            Given an \emph{isomorphism} $f\colon E\congto F$ of vector bundles on $X$, along with a connection $\nabla_F$ on $F$, we define the \define{pullback connection} $f^*\nabla_F$ on $E$ by
            \[
                f^*\nabla_F = (f^{-1}\otimes f^*)\circ\nabla_F\circ f
            \]
            (where $f^*$ denotes the pullback of differential forms).
            Locally, this takes the trivial connection $\d$ on $F$ to the connection $\d+f^{-1}\d f$ on $E$.

            Note that this is \emph{different} from the `other' notion of the pullback of a connection, which normally means pulling back along some change of base $f\colon X\to Y$.
        \end{definition}

        \begin{example}[Green's example (continued)]\label{example:green's-example-cont}
            The prototypical example of what we define in this section is Green's `barycentric connection', which can be understood as a way of building a connection $\nabla_p$ on $\overline{\mathcal{E}^p}$ given local connections $\nabla_{\alpha_i}$ on each $\mathcal{E}^0$, by simply defining $\nabla_p=\sum_{i=0}^p t_i\nabla_{\alpha_i}$.
            We explain this in slightly more detail by continuing \cref{example:green's-example}.

            \medskip

            Assuming that we have a basis of local sections over $U_1$, and another over $U_2$, we can look at all the connections we have on the $\mathcal{E}^{p,i}$ given by using the isomorphism $A_{12}^\bullet$ to pull back the local trivial connections.
            This data is given in \cref{table:pullback-connections-in-green's-example}.
            
            \begin{table}[ht]
                \centering
                \begin{tabular}{lllll}
                    bundle & $p$ & $p$-intersection & local connection &\\
                    \toprule
                    \multirow{4}{*}{$\mathcal{E}^{p,0}$} & \multirow{2}{*}{0} & $U_1$ & $\d$\\
                    & & $U_2$ & $\d$\\
                    & \multirow{2}{*}{1} & \multirow{2}{*}{$U_{12}$} & $\d$ & (from $U_1$)\\
                    & & & $\d+(A_{1,2}^0)^{-1}\d A_{1,2}^0=\d$ & (from $U_2$)\\
                    \midrule
                    \multirow{4}{*}{$\mathcal{E}^{p,1}$} & \multirow{2}{*}{0} & $U_1$ & $\d$\\
                    & & $U_2$ & $\d$\\
                    & \multirow{2}{*}{1} & \multirow{2}{*}{$U_{12}$} & $\d$ & (from $U_1$)\\
                    & & & $\d+(A_{1,2}^1)^{-1}\d A_{1,2}^1=\d+\frac{\d z}{z}$ & (from $U_2$)\\
                    \bottomrule
                \end{tabular}
                \caption{All the local connections for this example.}\label{table:pullback-connections-in-green's-example}
            \end{table}

            Using these local connections, we can form the \define{barycentric connections} $\nabla_\bullet^i$ on $\mathcal{E}^{\bullet,i}$ as follows:
            \begin{align*}
                \nabla_\bullet^0\text{ on }\mathcal{E}^{\bullet,0}\text{ is given by }
                &\begin{cases}
                    \nabla_0^0 = t_0\d&=\d\\
                    \nabla_1^0 = t_0\d+t_1\d&=\d
                \end{cases}\\
                \nabla_\bullet^1\text{ on }\mathcal{E}^{\bullet,1}\text{ is given by }
                &\begin{cases}
                    \nabla_0^1 = t_0\d&=\d\\
                    \nabla_1^1 = t_0\d+t_1\left(\d+\frac{\d z}{z}\right)&=\d+t_1\frac{\d z}{z}.
                \end{cases}
            \end{align*}

            \medskip

            We could continue this example, using \cref{example:fibre-integration-gives-cech-de-rham} to calculate an explicit Čech representative of the Chern class of $\mathscr{F}$ in de Rham cohomology, but this is not within the goals of this current paper; see instead the sequel to this paper.
        \end{example}

        \begin{definition}\label{definition:simplicial-connection}
            A \define{simplicial connection} $\nabla_\bullet$ on a vector bundle on the nerve $\mathcal{E}^\bullet$ is a family $\nabla_\bullet=\{\nabla_p\}_{p\in\mathbb{N}}$ of connections, where $\nabla_p$ is a connection on $\overline{\mathcal{E}^p}$, such that the comparison maps
            \[
                \left(\nerve{\bullet}f_p^i\times\id\right)^*
                \left(\overline{\mathcal{E}^{p-1}},\nabla_{p-1}\right)
                \xrightarrow{\comparison{p}^i(\mathcal{E}^\bullet)}
                \left(\id\times f_p^i\right)^*
                \left(\overline{\mathcal{E}^p},\nabla_p\right)
            \]
            are true morphisms.

            N.B. The pullbacks act on both the vector bundle and the connection \emph{simultaneously}: we do \emph{not} mean e.g. `pull back the vector bundle on the nerve by $\nerve{\bullet}f_p^i$ and the connection by $\id$'; we mean `pullback \emph{both} the vector bundle on the nerve and the connection by $(\nerve{\bullet}f_p^i\times\id)$'.
        \end{definition}

        \begin{remark}
            The individual connections $\nabla_p$ making up a simplicial connection $\nabla_\bullet$ on $\mathcal{E}^\bullet$ are connections on $\overline{\mathcal{E}^p}=\pi_p^*\mathcal{E}^p$, \emph{not} on $\mathcal{E}^p$ itself.
        \end{remark}

        \begin{remark}
            We ask in \cref{definition:simplicial-connection} that the morphisms between simplicial levels of the connection be true morphisms, but when it comes to defining the category of vector bundles on the nerve endowed with simplicial connections, a morphism $f\colon(\mathcal{E}^\bullet,\nabla_\bullet)\to(\mathcal{E}'^{\bullet},\nabla'_\bullet)$ between such objects will \emph{not} be asked to satisfy the corresponding requirement: it will simply be a morphism $f\colon\mathcal{E}^\bullet\to\mathcal{E}'^\bullet$ of vector bundles on the nerve.
        \end{remark}

        \begin{lemma}
            The curvature (defined simplicial degree by simplicial degree) of a simplicial connection is an endomorphism-valued simplicial 2-form.

            \begin{proof}
                This is a direct consequence of the definitions, as well as the aforementioned fact (in \cref{definition:true-morphism}) that a morphism that is true with respect to connections is also true with respect to their curvatures.
            \end{proof}
        \end{lemma}

        \begin{definition}\label{definition:admissible-simplicial-connection}
            We say that a simplicial connection is \define{admissible} if its curvature is an \emph{admissible} endomorphism-valued simplicial 2-form.
        \end{definition}

        \begin{lemma}\label{lemma:admissible-simplicial-connection-criterion}
            For a simplicial connection $\nabla_\bullet$ on $\mathcal{E}^\bullet$ to be admissible, it is sufficient to ask for sub-bundles $A^p,B^p\hookrightarrow\overline{\mathcal{E}^{p}}$ such that
            \begin{enumerate}[(i)]
                \item $A^p$ and $B^p$ are $\nabla_{p}$-flat;
                \item the comparison map
                    \begin{equation*}
                        (\nerve{\bullet} f_p^i\times\id)^*\left(\overline{\mathcal{E}^{p-1}},\nabla_{p-1}\right)
                        \xrightarrow{\comparison{p}^i}
                        (\id\times f_p^i)^*\left(\overline{\mathcal{E}^p},\nabla_p\right).
                    \end{equation*}
                    (which is already known to be a \emph{true} morphism, since the connection is assumed to be simplicial) restricts to a morphism
                    \begin{equation*}
                        (\nerve{\bullet} f_p^i\times\id)^*\left(A^{p-1},\nabla_{p-1}\right)
                        \xrightarrow{\comparison{p}^i}
                        (\id\times f_p^i)^*\left(B^p,\nabla_p\right);
                    \end{equation*}
                \item the above restriction of the comparison map induces an isomorphism
                \begin{equation*}
                    \comparison{p}^i\colon
                    (\nerve{\bullet} f_p^i\times\id)^*\left(\overline{\mathcal{E}^{p-1}}/A^{p-1}\right)
                    \congto
                    (\id\times f_p^i)^*\left(\overline{\mathcal{E}^p}/B^p\right).
                \end{equation*}
            \end{enumerate}

            \begin{proof}
                If $A^p$ (resp. $B^p$) is $\nabla_p$-flat then, in particular, it lies in the kernel of $\kappa(\nabla_p)$.
                Since $\kappa(\nabla_p)$ is simply $\nabla_p\circ\nabla_p$, the (again) aforementioned fact (that a morphism that is true with respect to some connections is also true with respect to their curvatures) tells us that the (true) morphism
                \begin{equation*}
                    (\nerve{\bullet} f_p^i\times\id)^*\left(A^{p-1},\nabla_{p-1}\right)\to(\id\times f_p^i)^*\left(B^p,\nabla_p\right)
                \end{equation*}
                induces a \emph{true} morphism
                \begin{equation*}
                    (\nerve{\bullet} f_p^i\times\id)^*\left(A^{p-1},\kappa(\nabla_{p-1})\right)\to(\id\times f_p^i)^*\left(B^p,\kappa(\nabla_p)\right).
                \end{equation*}
                But then, by \cref{lemma:admissible-endomorphism-valued-form-criterion}, we are done.
            \end{proof}
        \end{lemma}
        
        \begin{definition}\label{definition:compatible-family-of-simplicial-connections}
            The difference $\nabla'_\bullet-\nabla_\bullet$ of two admissible simplicial connections has \emph{no a priori reason} to be an admissible endomorphism-valued simplicial 1-form, which prompts the following definition: a set of admissible simplicial connections is said to be \define{compatible} if the difference of any two simplicial connections is indeed an admissible endomorphism-valued simplicial 1-form.
        \end{definition}
        
        \begin{remark}
            In summary of the definitions in this chapter so far:
            \begin{itemize}
                \item the \emph{simplicial} condition (\cref{definition:simplicial-connection}) ensures that various forms (such as the curvature) defined by a connection satisfy the gluing condition needed to give a simplicial differential form;
                \item the \emph{admissibility} condition (\cref{definition:admissible-simplicial-connection}) will ensure that we can evaluate `generalised invariant polynomials' on the curvature and get something that is the same in \emph{all} simplicial degrees;
                \item the \emph{compatibility} condition (\cref{definition:compatible-family-of-simplicial-connections}) will ensure that characteristic classes will be independent of the choice of connection.
            \end{itemize}
            The last two points will be further explained and justified in the sequel to this paper.
        \end{remark}

    \subsection{Being generated in degree zero}

        \begin{remark}\label{remark:explaining-generated-in-degree-zero-definition}
            For an arbitrary vector bundle on the nerve $\mathcal{E}^\bullet$, there is no reason for $\mathcal{E}^0$ and $\mathcal{E}^p$ to be isomorphic.
            We might imagine, however, that there should be some sort of structural condition ensuring that $\mathcal{E}^0_\alpha$ and $\mathcal{E}^p_{\alpha\ldots\alpha}$ agree, so that the vector bundle is somehow `built up' from its degree-zero part.
            Indeed, Green bundles can be thought of exactly in this way: the cokernel of any $\mathcal{E}^\bullet f_p^i$ is elementary, and these cokernels satisfy a cocycle condition.
            If we write $\zeta_p^i\colon[0]\to[p]$ to mean the morphism in $\Delta$ that sends $0$ to $i$, then we have an \emph{equality} of sheaves
            \begin{equation}
            \label{equation:equality-of-sheaves-pullback-by-zeta-double-restriction}
                \left(\nerve{\bullet}\zeta_p^i\right)^*
                \mathcal{E}^0
                \restricted U_{\alpha_0\ldots\alpha_p}
                =
                \left(\mathcal{E}^0 \restricted U_{\alpha_i}\right)
                \restricted U_{\alpha_0\ldots\alpha_p}
            \end{equation}
            as well as an isomorphism (\cref{theorem:green's-resolution}~{(iv)})
            \begin{equation}
            \label{equation:splitting-of-Ep-before-pi-p}
                \mathcal{E}^p \restricted U_{\alpha_0\ldots\alpha_p}
                \cong
                \left( \left(\nerve{\bullet}\zeta_p^i\right)^* \mathcal{E}^0 \restricted U_{\alpha_0\ldots\alpha_p}\right)
                \oplus
                \mathcal{K}_{0,i,p}
            \end{equation}
            where $\mathcal{K}_{0,i,p} = \Coker(\mathcal{E}^\bullet\zeta_p^i)$ is elementary in the $\mathcal{E}_{\alpha_i}^0$.
            This splitting lets us think of $\mathcal{E}^p$ as being built up from $\mathcal{E}^0$ by adding an elementary sequence, and the cocycle condition (\cref{theorem:green's-resolution}~{(v)}) satisfied by the $\mathcal{K}_{0,i,p}$ tells us that building $\mathcal{E}^p$ from $\mathcal{E}^0$ and then $\mathcal{E}^{q}$ from $\mathcal{E}^p$ (for $q>p$) is the same as building $\mathcal{E}^q$ from $\mathcal{E}^0$ directly.
        \end{remark}

        \begin{remark}\label{remark:explaining-compatible-extension-of-connection}
            Applying $\pi_p^*$ (which is exact) to \cref{equation:splitting-of-Ep-before-pi-p}, we get a splitting
            \begin{equation}
            \label{equation:that-splitting-of-Ep}
                \pi_p^* \mathcal{E}^p \restricted U_{\alpha_0\ldots\alpha_p}
                \cong
                \pi_p^* \left( \left(\nerve{\bullet}\zeta_p^i\right)^* \mathcal{E}^0 \restricted U_{\alpha_0\ldots\alpha_p}\right)
                \oplus
                \pi_p^* \mathcal{K}_{0,i,p}.
            \end{equation}
            So if we have local connections $\widetilde{\nabla}_{\alpha_i}$ on $\mathcal{E}^0\restricted U_{\alpha_i}$ for $i=0,\ldots,p$, then \cref{equation:equality-of-sheaves-pullback-by-zeta-double-restriction} tells us that $\pi_p^*\widetilde{\nabla}_{\alpha_i}$ defines a connection on the first summand.
            But we can also construct a connection on  the elementary summand by taking a direct sum of the $\widetilde{\nabla}_{\alpha_i}$.
            Together then, this lets us build a connection on each $\pi_p^*\mathcal{E}^p$.

            Importantly, this direct sum of connections on the elementary summand gives a \define{compatible\footnote{Not to be confused with a compatible \emph{family} of \emph{simplicial} connections.} sequence of connections}: the connections commute with the differentials.
            This is because being an elementary sequence means that all differentials are identity maps, so if we place the same connection in the two non-zero degrees of each elementary component (i.e. on both copies of $M$ in $0\to M\to M[1]\to 0$) then commutativity is trivial.
            But it is known, by applying \cite[Lemma~4.22]{Baum&Bott1972} to the exact sequence $0\to M \to (M\to M[1])\to M[1]\to0$, that such a sequence of connections gives a trivial characteristic class, and so extending by such a thing will ensure that we do not change the characteristic class of our bundle, thanks to additivity.
        \end{remark}

        \begin{definition}\label{definition:generated-in-degree-zero}
            A simplicial connection $\nabla_\bullet$ on a \emph{Gre{}en} vector bundle on the nerve $\mathcal{E}^\bullet$ is said to be \define{generated in degree zero} if, for all $U_{\alpha_0\ldots\alpha_p}\in\cover$, it is of the form
            \[
                \nabla_p\restricted U_{\alpha_0\ldots\alpha_p}
                = \sum_{i=0}^p t_i \pi_p^*\widetilde{\nabla}_{\alpha_i}
            \]
            where each $\widetilde{\nabla}_{\alpha_i}$ is a connection on $\mathcal{E}^0\restricted U_{\alpha_i}$.
            That is, if it is built as in \cref{remark:explaining-compatible-extension-of-connection}.
        \end{definition}

        \begin{remark}
            In \cite{Hosgood2020}, we try to use the terminology `barycentric' only to describe connections generated in degree zero on pullbacks (to the nerve) of global vector bundles; the phrase `generated in degree zero' applies to arbitrary vector bundles on the nerve.
            A mild abuse of language is permitted, however, by the fact that, if the bundle $\mathcal{E}^\bullet$ is strongly cartesian (which is the case for pullbacks of global bundles), then the inclusion maps $\mathcal{E}^\bullet\zeta_p^i$ are isomorphisms.
            This means, in particular, that the definition of ‘being generated in degree zero’ agrees with that of the barycentric connection on pullbacks of global vector bundles.
        \end{remark}

        \begin{theorem}\label{theorem:admissible-gidz-on-Green}
            Let $\mathcal{E}^{\bullet,\anotherbullet}$ be a Gre{}en complex.
            Then we can endow each $\mathcal{E}^{\bullet,j}$ with a simplicial connection generated in degree zero.
            Further, these simplicial connections are admissible.

            \begin{proof}
                We split the proof into three steps: \emph{defining} connections that are \emph{generated in degree zero}; showing that they are \emph{simplicial}; and then showing that they are \emph{admissible}.
                For ease of notation, we write $\mathcal{E}^\bullet$ instead of $\mathcal{E}^{\bullet,j}$.
                \begin{enumerate}
                    \item Take \emph{arbitrary} local connections $\widetilde{\nabla}_\alpha$ on $\mathcal{E}^0\restricted U_\alpha$ for all $\alpha$.
                        Since $\mathcal{E}^{\bullet,\anotherbullet}$ is Gre{}en, we can use \cref{remark:explaining-generated-in-degree-zero-definition,remark:explaining-compatible-extension-of-connection} to define connections $\nabla_\bullet$ on $\overline{\mathcal{E}^\bullet}$ that are \emph{generated in degree zero}.
                        \label{item:we-can-define-gen-in-deg-zero-connections}
                    \item To show that each $\nabla_\bullet$ defined above is a \emph{simplicial} connection, we need to show that
                        \begin{equation}
                        \label{equation:simplicial-condition-square}
                            \begin{tikzcd}[column sep=huge,row sep=7em]
                                \left(\nerve{\bullet} f_p^i\times\id\right)^*\pi_{p-1}^*{\mathcal{E}^{p-1}}
                                    \ar[r,hook,"\comparison{p}^i(\mathcal{E}^\bullet)"]
                                    \ar[d,swap,"\left(\nerve{\bullet} f_p^i\times\id\right)^*\nabla_{p-1}"]
                                &\left(\id\times f_p^i\right)^*\pi_p^*{\mathcal{E}^p}
                                    \ar[d,"\left(\id\times f_p^i\right)^*\nabla_p"]\\
                                \begin{tabular}{c}
                                    $\left(\nerve{\bullet} f_p^i\times\id\right)^*\pi_{p-1}^*{\mathcal{E}^{p-1}}$\\
                                    $\otimes\,\Omega^1_{\nervesimplex{p}}$
                                \end{tabular}
                                    \ar[r,hook,swap,"\comparison{p}^i(\mathcal{E}^\bullet)\,\otimes\,\id"]
                                &\begin{tabular}{c}
                                    $\left(\id\times f_p^i\right)^*\pi_p^*{\mathcal{E}^p}$\\
                                    $\otimes\,\Omega^1_{\nervesimplex{p}}$
                                \end{tabular}
                            \end{tikzcd}
                        \end{equation}
                        commutes (where we use \cref{corollary:Green-has-injective-comparisons}, which tells us that the comparison maps are injective, as well as the fact that tensoring preserves splittings, to see that the horizontal arrows are injections).
                        We start by making some simplifications.

                        First of all, all of the sheaves in \cref{equation:simplicial-condition-square} lie over $\nerve{p}\times\Delta^{p-1}$, but we make the identification
                        \begin{align*}
                            \nerve{p}\times\Delta^{p-1}
                            &\simeq
                            \nerve{p}\times f_p^i\left(\Delta^{p-1}\right)\\
                            &\subset
                            \nervesimplex{p}
                        \end{align*}
                        so that we can label \emph{both} (i.e. on the nerve and on the simplex) simplicial parts with the same indices: the nerve being labelled with $\{0,1,\ldots,p\}$; the simplex being labelled with $\{0,1\ldots,\hat{i},\ldots,p\}$ (so that $i$ is now fixed).

                        Next, we can use the commutativity of the square
                        \begin{equation}
                        \label{equation:square-for-pullback-restriction-equality}
                            \begin{tikzcd}[row sep=huge]
                                \nerve{p}
                                    \ar[r,"\nerve{\bullet}f_p^i"]
                                & \nerve{p-1}\\
                                U_{\alpha_0\ldots\alpha_p}
                                    \ar[r,hook]
                                    \ar[u,hook]
                                & U_{\alpha_0\ldots\widehat{\alpha_i}\ldots\alpha_p}
                                    \ar[u,hook]
                            \end{tikzcd}
                        \end{equation}
                        to see that
                        \begin{align*}
                            \left(\left(\nerve{\bullet}f_p^i\times\id\right)^* \overline{\mathcal{E}^{p-1}}\right) &\restricted \left(U_{\alpha_0\ldots\alpha_p}\times f_p^i(\Delta^{p-1})\right)\\
                            =
                            \left(\overline{\mathcal{E}^{p-1}} \restricted \left(U_{\alpha_0\ldots\widehat{\alpha_i}\ldots\alpha_p}\times\Delta^{p-1}\right)\right) &\restricted \left(U_{\alpha_0\ldots\alpha_p}\times f_p^i(\Delta^{p-1})\right)
                        \end{align*}
                        as sheaves over $U_{\alpha_0\ldots\alpha_p}\times f_p^i(\Delta^{p-1})$.

                        So, to prove the commutativity of \cref{equation:simplicial-condition-square}, we start by calculating how the two pullbacks of the connections (that is, the vertical arrows in the square) act (after restricting all sheaves to $U_{\alpha_0\ldots\alpha_p}$).
                        The first pullback is
                        \begin{equation*}
                            \left(\nerve{\bullet}f_p^i\times\id\right)^*\nabla_{p-1}
                            =
                            \left(\nerve{\bullet}f_p^i\times\id\right)^* \sum_{\substack{j=0\\j\neq i}}^p t_j\pi_{p-1}^*\widetilde{\nabla}_{\alpha_j}
                            =
                            \sum_{\substack{j=0\\j\neq i}}^p t_j\pi_{p}^*\widetilde{\nabla}_{\alpha_j}
                        \end{equation*}
                        where the first equality is simply the definition, and the second equality is the more subtle one.
                        We are really writing $\widetilde{\nabla}_{\alpha_j}$ to mean two different things: on the left-hand side, it means the connection $\pi_{p-1}^*\widetilde{\nabla}_{\alpha_j}$ extended by the trivial connection (as explained in \cref{remark:explaining-compatible-extension-of-connection}) on $\mathcal{K}_{0,j,p-1}$; on the right-hand side, it means the connection $\pi_p^*\widetilde{\nabla}_{\alpha_j}$ extended by the trivial connection on $\mathcal{K}_{0,j,p}$ (all of this using the notation and properties of \cref{remark:explaining-generated-in-degree-zero-definition,remark:explaining-compatible-extension-of-connection}).
                        So this second equality will really follow from \cref{equation:cocycle-condition-on-the-K}, which tells us that extending trivially on $\mathcal{K}_{0,j,p-1}$ and then again to the rest of $\mathcal{K}_{0,j,p}$ (which we do for $\pi_{p-1}^*\widetilde{\nabla}_{\alpha_j}$) is the same as simply extending trivially on $\mathcal{K}_{0,j,p}$ (which we do for $\pi_{p}^*\widetilde{\nabla}_{\alpha_j}$).

                        The second pullback is much simpler:
                        \begin{equation*}
                            \left(\id\times f_p^i\right)^* \nabla_p
                            =
                            \left(\id\times f_p^i\right)^* \sum_{j=0}^p t_j\pi_p^*\widetilde{\nabla}_{\alpha_j}
                            =
                            \sum_{\substack{j=0\\j\neq i}}^p t_j\pi_p^*\widetilde{\nabla}_{\alpha_j}
                        \end{equation*}
                        which is `exactly the same' as the first pullback --- the scare quotes being important, because these connections are on \emph{different sheaves}.
                        But, since the horizontal arrows in \cref{equation:simplicial-condition-square} are injections, it means that these two connections really are the same when we just follow how they act on the top-left sheaf in the square; i.e., the square commutes.

                        \medskip

                        Looking ahead to \cref{equation:Ep-decomposed-into-Ep-1}, since we have extended by something \emph{compatible} on $\mathcal{K}_{p-1,i,p}$, the characteristic class of $\mathcal{K}_{p-1,i,p}$ will be 1, as mentioned above.
                        This means, by additivity of characteristic classes, that the classes of $\mathcal{E}^{p-1}\restricted U_{\alpha_0\ldots\alpha_p}$ and $\mathcal{E}^p\restricted U_{\alpha_0\ldots\alpha_p}$ will agree.
                        This is the content of half of the proof of \cite[Lemma~2.2]{Green1980}.
                    \item To show that each $\nabla_\bullet$ is an \emph{admissible} simplicial connection, it suffices to show that the conditions in \cref{lemma:admissible-simplicial-connection-criterion} are satisfied.

                        Before proceeding with the proof, however, we take some time to look at how the splittings that we have been using respect the simplicial structure.
                        The Gre{}en assumption implies that
                        \begin{equation}
                        \label{equation:Ep-decomposed-into-Ep-1}
                            \mathcal{E}^p \restricted U_{\alpha_0\ldots\alpha_p}
                            \cong
                            \left(
                                \mathcal{E}^{p-1} \restricted U_{\alpha_0\ldots\widehat{\alpha_i}\ldots\alpha_p}
                            \right) \restricted U_{\alpha_0\ldots\alpha_p}
                            \oplus
                            \mathcal{K}_{p-1,i,p}
                        \end{equation}
                        which, combined with \cref{equation:that-splitting-of-Ep}, says that (for $i\neq j$)\footnote{We can ignore the case where $i=j$, since we will always be working on the embedding of $\Delta^{p-1}$ into $\Delta^p$ by $f_p^i$, which is equivalent to working with the labelling $(\alpha_0,\ldots,\widehat{\alpha_i},\ldots,\alpha_p)$ on $\Delta^{p-1}$.}
                        \begin{align*}
                            \mathcal{E}^p \restricted U_{\alpha_0\ldots\alpha_p}
                            &\cong
                            \left(
                                \big(\nerve{\bullet}\zeta_{p-1}^j\big)^* \mathcal{E}^0 \restricted U_{\alpha_0\ldots\widehat{\alpha_i}\ldots\alpha_p}
                                \oplus
                                \mathcal{K}_{0,j,p-1}
                            \right) \restricted U_{\alpha_0\ldots\alpha_p}
                            \oplus
                            \mathcal{K}_{p-1,i,p}\\
                            &\cong
                            \left(
                                \big(\nerve{\bullet}\zeta_{p-1}^j\big)^* \mathcal{E}^0 \restricted U_{\alpha_0\ldots\widehat{\alpha_i}\ldots\alpha_p}
                            \right) \restricted U_{\alpha_0\ldots\alpha_p}
                            \oplus
                            \left(
                                \mathcal{K}_{0,j,p-1} \restricted U_{\alpha_0\ldots\alpha_p}
                                \oplus
                                \mathcal{K}_{p-1,i,p}
                            \right)
                        \end{align*}
                        and (again by the Gre{}en assumption) these $\mathcal{K}$ satisfy some cocycle condition:
                        \begin{equation}
                        \label{equation:cocycle-condition-on-the-K}
                            \mathcal{K}_{0,i,p}
                            \cong
                            \mathcal{K}_{0,j,p-1} \restricted U_{\alpha_0\ldots\alpha_p}
                            \oplus
                            \mathcal{K}_{p-1,i,p}
                        \end{equation}
                        whence
                        \begin{equation}
                        \label{equation:Ep-decomposed-into-Ep-1-as-E0-pullback}
                            \mathcal{E}^p \restricted U_{\alpha_0\ldots\alpha_p}
                            \cong
                            \left(
                                \big(\nerve{\bullet}\zeta_{p-1}^j\big)^* \mathcal{E}^0 \restricted U_{\alpha_0\ldots\widehat{\alpha_i}\ldots\alpha_p}
                            \right) \restricted U_{\alpha_0\ldots\alpha_p}
                            \oplus
                            \mathcal{K}_{0,i,p}.
                        \end{equation}

                        Now, for each comparison map $\comparison{p}^i$, we set
                        \begin{align*}
                            A^{p-1} &= \pi_p^*\mathcal{K}_{0,j,p-1}\\
                            B^p &= \pi_p^*\mathcal{K}_{0,i,p}
                        \end{align*}
                        for an arbitrary $j\neq i$.
                        \begin{enumerate}[(i)]
                            \item $A^p$ is $\nabla_p$ flat by definition: we extended the connection $\pi_p^*\widetilde{\nabla}_{\alpha_i}$ by the trivial connection on this direct summand (cf. above, or \cref{remark:explaining-compatible-extension-of-connection}).
                            \item The comparison map $\comparison{p}^i$ is simply the pullback (along $\pi_{p,p-1}$) of $\mathcal{E}^\bullet f_p^i$, and so respects the splitting \cref{equation:that-splitting-of-Ep} almost by definition, because these splittings come from the Gre{}en assumption.
                            \item We want the comparison map to induce an isomorphism when we take the quotients.
                                To avoid getting lost in a mire of notation, instead of using \cref{equation:that-splitting-of-Ep}, we use the language of \cref{definition:Gre{}en}.
                                Set
                                \begin{align*}
                                    \upalpha &= (\alpha_0,\ldots,\alpha_p)\\
                                    \upbeta &= (\alpha_0,\ldots,\widehat{\alpha_i},\ldots,\alpha_p)\\
                                    \upgamma &= (\alpha_j)
                                \end{align*}
                                and note that $A^{p-1}=\pi_{p-1}^* \mathcal{L}_{\upbeta,\upgamma}$ and $B^p=\pi_p^* \mathcal{L}_{\upalpha,\upgamma}$.
                                Then
                                \begin{equation*}
                                    \begin{array}{lr}
                                        \mathcal{E}_\upalpha^\bullet
                                        \cong
                                        \mathcal{E}_\upgamma^\bullet \restricted U_\upbeta \restricted U_\upalpha
                                        &\oplus
                                        \mathcal{L}_{\upalpha,\upgamma}\\[0.3em]
                                        \mathcal{E}_\upbeta^\bullet
                                        \cong
                                        \mathcal{E}_\upgamma^\bullet \restricted U_\upbeta
                                        &\oplus
                                        \mathcal{L}_{\upbeta,\upgamma}
                                    \end{array}
                                \end{equation*}
                                and we are interested in
                                \begin{equation*}
                                    \left(\nerve{\bullet}f_p^i\times\id\right)^* \pi_{p-1}^* \left(\mathcal{E}_\upgamma^\bullet \restricted U_\upbeta\right)
                                    \xrightarrow{\pi_{p,p-1}^*\left(\mathcal{E}^\bullet f_p^i\right)}
                                    \left(\id\times f_p^i\right)^* \pi_p^* \left(\mathcal{E}_\upgamma^\bullet \restricted U_\upbeta \restricted U_\upalpha\right).
                                \end{equation*}

                                Consider first the source of this map: by \cref{definition:comparison-map}, we know that
                                \begin{equation*}
                                    \left(\nerve{\bullet}f_p^i\times\id\right)^* \pi_{p-1}^* \left(\mathcal{E}_\upgamma^\bullet \restricted U_\upbeta\right)
                                    \cong
                                    \pi_{p,p-1}^* \left(\nerve{\bullet}f_p^i\right)^* \left(\mathcal{E}_\upgamma^\bullet \restricted U_\upbeta\right)
                                \end{equation*}
                                but, as in \cref{equation:square-for-pullback-restriction-equality}, we know that
                                \begin{equation*}
                                    \left(\nerve{\bullet}f_p^i\right)^* \left(\mathcal{E}_\upgamma^\bullet \restricted U_\upbeta\right) \restricted U_\upalpha
                                    =
                                    \left(\mathcal{E}_\upgamma^\bullet \restricted U_\upbeta\right) \restricted U_\upbeta \restricted U_\upalpha
                                    =
                                    \mathcal{E}_\upgamma^\bullet \restricted U_\upbeta \restricted U_\upalpha
                                \end{equation*}

                                The target is much simpler, since, as explained in \cref{definition:comparison-map}, the pullback along $f_p^i$ on the simplex part changes nothing:
                                \begin{equation*}
                                    \left(\id\times f_p^i\right)^* \pi_p^* \left(\mathcal{E}_\upgamma^\bullet \restricted U_\upbeta \restricted U_\upalpha\right)
                                    \cong
                                    \pi_{p,p-1}^* \left(\mathcal{E}_\upgamma^\bullet \restricted U_\upbeta \restricted U_\upalpha\right).
                                \end{equation*}
                                But then we see that the source and target, \emph{when we restrict to $U_\alpha$}, are identical, so it would suffice to show that the comparison map descends to an injection when we take these quotients.
                                But these quotients are exactly the sheaf with which we started: we added the cokernels and then quotiented them out; and we know, from \cref{corollary:Green-has-injective-comparisons}, that the comparison map is injective here.
                        \end{enumerate}
                \end{enumerate}
            \end{proof}
        \end{theorem}

        \begin{remark}
            We can, in particular, think of \cref{theorem:admissible-gidz-on-Green} as a proof that, \emph{when working with Gre{}en complexes}, being generated in degree zero implies admissibility.
        \end{remark}

        \begin{lemma}\label{lemma:green-gidz-compatible-family}
            Let $\mathcal{E}^\bullet$ be a Gre{}en vector bundle on the nerve, and let $\nabla_\bullet$ and $\nabla'_\bullet$ be two simplicial connections on $\mathcal{E}^\bullet$ that are generated in degree zero.
            Then the difference $\nabla'_\bullet-\nabla_\bullet$ is an admissible endomorphism-valued simplicial form.
            In other words, the set of generated-in-degree-zero connections on a Green vector bundle on the nerve is a compatible family.

            \begin{proof}
                The fact that the difference of two arbitrary simplicial connections on an arbitrary vector bundle on the nerve is an endomorphism-valued simplicial 1-form follows ``immediately'' from the definitions, without any extra hypotheses.
                The content of this lemma is that \emph{generated-in-degree-zero} connections on a \emph{Green} vector bundle on the nerve have an \emph{admissible} difference.

                Write $\nabla_\bullet^{(1)}$ and $\nabla_\bullet^{(2)}$ to mean the two connections, so that $\nabla_\bullet^{(i)}$ is constructed as in \cref{remark:explaining-generated-in-degree-zero-definition}, but with a different choice of local connections $\widetilde{\nabla}_\alpha^{(i)}$ on each $\mathcal{E}_0\restricted U_\alpha$.
                This means that, for each $p\in\mathbb{N}$, we can write
                \begin{equation*}
                    \nabla_p^{(i)}
                    = \sum_{j=0}^p t_j\pi_p^*\widetilde{\nabla}^{(i)}_{\alpha_j}
                \end{equation*}
                which means that the difference is of the form
                \begin{equation*}
                    \nabla_p^{(2)} - \nabla_p^{(1)}
                    = \sum_{j=0}^p t_j \pi_p^* \left(\widetilde{\nabla}^{(2)}_{\alpha_j} - \widetilde{\nabla}^{(1)}_{\alpha_j}\right)
                    = \sum_{j=0}^p t_j \pi_p^* \eta_{\alpha_j}
                    = \sum_{j=0}^p t_j \overline{\eta_{\alpha_j}}
                \end{equation*}
                where $\eta_{\alpha_j}$ is an endomorphism-valued simplicial 1-form on $U_{\alpha_j}$.
                The claim is then that $\sum_{j=0}^p t_j \overline{\eta_{\alpha_j}}$ is an \emph{admissible} endomorphism-valued simplicial 1-form, and this follows almost exactly as in the proof of \cref{theorem:admissible-gidz-on-Green}, because the difference $\sum t_j\overline{\eta_{\alpha_j}}$ is somehow also ``generated in degree zero'', in that it is given by trivially extending on each $(\nerve{\bullet}f_p^i)^*\mathcal{E}_{p-1}\hookrightarrow\mathcal{E}_p$.
            \end{proof}
        \end{lemma}

    \subsection{Equivalences}

        Now we can improve upon the result of \cref{corollary:equivalence-without-connections}.

        \begin{definition}\label{definition:cartvectzeroX-and-greenzeroX}
            We define the category $\greenzeroX$ via the Grothendieck construction applied to the functor $F\colon\greenX\to\mathsf{Set}$ given, on an object $\mathcal{E}^{\bullet,\anotherbullet}$, by
            \[
                F(\mathcal{E}^{\bullet,j}) = \big\{ \mbox{generated-in-degree-zero simplicial connections on $\mathcal{E}^{\bullet,j}$} \big\},
            \]
            where \cref{theorem:admissible-gidz-on-Green} tells us that this set is non-empty.
            So an object of $\greenzeroX$ is a pair $(\mathcal{E}^{\bullet,\anotherbullet},\nabla_\bullet^\anotherbullet)$, where $\mathcal{E}^{\bullet,\anotherbullet}$ is an object of $\greenX$, and $\nabla_\bullet^j$ is a simplicial connection generated in degree zero (and thus admissible) on $\mathcal{E}^{\bullet,j}$; the morphisms of $\greenzeroX$ are exactly those of $\greenX$.
            In particular, $\greenzeroX$ is a homotopical category with the same weak equivalences as $\greenX$
        \end{definition}

        \begin{remark}\label{remark:4-is-an-equivalence-of-1-cats}
            By construction, the forgetful functor
            \[
                \numberincircle{4}\colon
                \greenzeroX \to \greenX
            \]
            (that forgets about the connections) is fully faithful and essentially surjective, and thus induces an equivalence of categories.
            More importantly, though, it also restricts to give an equivalence of the corresponding subcategories of weak equivalences, which will be useful later on.
        \end{remark}

        \begin{remark}\label{remark:diagram-descends-to-infty-1-and-4-is-equivalence}
            Note that $\numberincircle{4}$ is constructed in such a way that it automatically preserves weak equivalences, as explained in \cref{definition:cartvectzeroX-and-greenzeroX}.

            In fact, by \cref{remark:4-is-an-equivalence-of-1-cats}, $\numberincircle{4}$ actually directly induces an equivalence at the level of localisations:
            \[
                \numberincircle{4}\colon
                \LL{\greenzeroX}
                \congto
                \LL{\greenX}.
            \]
            This is because an equivalence of relative categories that restricts to an equivalence of the wide subcategories of weak equivalences induces an equivalence of the localisations, as can be shown by using \cite[Lemma~5.4]{Barwick&Kan2012}, as well as the fact that an equivalence of categories induces an equivalence of their nerves.
    
            As a side note, even though $\numberincircle{4}$ gives an equivalence at the level of localisations, it there no longer looks like a Grothendieck construction: there is no reason for two weakly equivalent Gre{}en complexes (or even, more simply, quasi-isomorphic complexes of vector bundles) to admit local connections that are in bijective correspondence with one another.
        \end{remark}

        \begin{corollary}\label{theorem:the-main-coherent-theorem}
            There is an equivalence of $(\infty,1)$-categories
            \begin{equation*}
                \hocolim_\cover\LL{\greenzeroX} \simeq \hocolim_\cover\LL{\gcohUX}.
            \end{equation*}
            \begin{proof}
                \cref{remark:diagram-descends-to-infty-1-and-4-is-equivalence} tells us that $\numberincircle{4}$ is an equivalence; \cref{corollary:equivalence-without-connections} tells us that $\numberincircle{3}\numberincircle{2}\numberincircle{1}$ is an equivalence.
            \end{proof}
        \end{corollary}

        \begin{lemma}\label{lemma:gcohUX-to-gccohX-essentially-surjective}
            The functor
            \[
                \hocolim_\cover\LL{\gcohUX}\to\LL{\gccohX}
            \]
            induced by the full embedding $\gcohUX\hookrightarrow\gccohX$ is essentially surjective.

            \begin{proof}
                We need to show that, given any $\mathscr{K}^\bullet$ with coherent cohomology, there exists some cover $\cover$ such that, on each $U_\alpha\in\cover$, there exists some bounded complex of coherent sheaves quasi-isomorphic to $\mathscr{K}^\bullet\restricted U_\alpha$.
                Following the proof of \cite[Proposition~1.7.11]{Kashiwara&Schapira1990}, we see that it suffices to show that, for any surjective morphism $u\colon\mathscr{G}\to\mathscr{H}$ of analytic sheaves with $\mathscr{H}$ coherent, and any point $x\in X$, there exists a neighbourhood $U_x$ of $x$, a coherent sheaf $\mathscr{F}$ on $U_x$, and a morphism $t\colon\mathscr{F}\to\mathscr{G}$ such that the composite morphism $ut\colon\mathscr{F}\to\mathscr{H}$ is surjective.

                To see this, let $x\in X$.
                Since $\mathscr{H}$ is coherent, $\mathscr{H}_x$ is of finite type over the local ring $\OO_{X,x}$, and so can be generated by a finite number of sections $s_1,\ldots,s_r$.
                Since the map $\mathscr{G}_x\to\mathscr{H}_x$ is surjective, we can lift these sections to sections $t_1,\ldots,t_r$ of $\mathscr{G}_x$ defined on some neighbourhood $U_x$ of $x$.
                Define the sheaf $\mathscr{F}$ on $U_x$ to be the free sheaf $(\OO_{U_x})^r$, and define the morphism $t\colon\mathscr{F}\to\mathscr{G}$ by the sections $t_1,\ldots,t_r$.
                Consider the cokernel $\mathscr{C}$ of $ut\colon\mathscr{F}\to\mathscr{H}$, which is coherent, since both $\mathscr{F}$ and $\mathscr{H}$ are.
                Since $\mathscr{C}$ is coherent, its support is an analytic set, and thus closed.
                By construction, the map $(ut)_x\colon\mathscr{F}\to\mathscr{H}$ is surjective.
                Since the complement of $\operatorname{supp}(\mathscr{C})$ is open, and not equal to $\{x\}$, there exists a neighbourhood $V_x$ of $x$ that does not intersect with $\operatorname{supp}(\mathscr{C})$.
                Thus $ut\colon\mathscr{F}\to\mathscr{H}$ is surjective on $V_x$.
            \end{proof}
        \end{lemma}

        \begin{lemma}\label{lemma:gcohUX-to-gccohX-fully-faithful}
            The functor
            \[
                \hocolim_\cover\LL{\gcohUX}\to\LL{\gccohX}
            \]
            induced by the full embedding $\gcohUX\hookrightarrow\gccohX$ is fully faithful.

            \begin{proof}
                If we identify each $\LL{\gcohUX}$ with its weakly-essential image in $\LL{\gccohX}$ (i.e. the subcategory spanned by objects \emph{weakly equivalent} to those in the image of the inclusion), then we can simply apply \cref{lemma:filtered-poset-colimit}.
            \end{proof}
        \end{lemma}

        \begin{definition}
            Let $\scartvectX$ be the subcategory of $\cartvectX$ spanned by complexes that actually \emph{are} cartesian complexes of locally free sheaves on the nerve; let $\sgreenzeroX$ be the subcategory of $\greenzeroX$ spanned by complexes that actually \emph{are} Green complexes.
        \end{definition}

        \begin{lemma}\label{lemma:can-sort-of-strictify-vect-and-green}
            $\hocolim_\cover\LL{\greenzeroX} \simeq \hocolim_\cover\LL{\sgreenzeroX}$.
            
            \begin{proof}
                The homotopy colimit over refinements of all covers is equivalent to the homotopy colimit over some truncation below of refinements over all covers, i.e. we can always assume that our covers are as fine as we wish when computing the homotopy colimit.
                But, by taking a fine enough cover $\cover$, every object of $\greenzeroX$ restricted to each $U_\alpha$ is free, and so, as in the proof of \cref{lemma:green-gives-essential-surjectivity}, we can invert quasi-isomorphisms whose target is in $\greenzeroX$.
                This means that, given some morphism in the localisation $\LL{\greenzeroX}$, expressed as a chain of roofs with the left-legs all quasi-isomorphisms, we can invert the quasi-isomorphisms and compose the resulting morphisms to obtain a single morphism in $\greenzeroX$ which is equal to that in the localisation with which we started.
                This means that $\LL{\greenzeroX}$ is equivalent to $\LL{\sgreenzeroX}$, and so their homotopy colimits agree.
            \end{proof}
        \end{lemma}

        \begin{corollary}\label{corollary:the-main-coherent-corollary}
            There is an equivalence of $(\infty,1)$-categories
            \begin{equation*}
                \hocolim_\cover\LL{\sgreenzeroX} \simeq \LL{\gccohX}.
            \end{equation*}
            \begin{proof}
                This is a consequence of \cref{theorem:the-main-coherent-theorem,lemma:gcohUX-to-gccohX-essentially-surjective,lemma:gcohUX-to-gccohX-fully-faithful,lemma:can-sort-of-strictify-vect-and-green}.
            \end{proof}
        \end{corollary}

\addcontentsline{toc}{section}{References}
\printbibliography

\end{document}